[1994/12/01]% LaTeX date must December 1994 or later

\documentclass[coversheets,math]{puthesis}

\usepackage{amscd}
\usepackage{amssymb}
\usepackage{amsmath}
\usepackage{amsfonts}
\usepackage{color}
\usepackage{graphicx}
\usepackage{verbatim}

\title{Error Correcting Codes on Algebraic Surfaces}

\author{Chris C. Lomont}{Lomont, Chris C.}

\degree{Doctor of Philosophy}{Ph.D.}{May}{2003}

\majorprof{Tzuong-Tsieng Moh}

\newcommand{\field}[1]{\ensuremath{\mathbb{F}_{#1}}} % finite field
\newcommand{\sheaf}[1]{\ensuremath{\mathcal{#1}}}    % sheaf
\newcommand{\proj}{\ensuremath{\mathbb{P}}}          % projective space
\newcommand{\pts}{\ensuremath{\gamma}}               % count of points on a curve
\newcommand{\intpts}{\ensuremath{\delta}}            % pts on intersection L.C_i
\newcommand{\Pic}{\ensuremath{\textrm{Pic}}}         % Picard group upright font

\newcommand{\blankpage}{{\clearpage {~} \vspace{4.0in} \begin{center}(This page deliberately left
blank.)\end{center}}}

\newcommand{\todo}[1]                                % change this to output all todo's

\newcommand{\qtn}[2]{{\small{\emph{``#2"} #1}}} % how to format quotations

\begin{document} % and we're off!!

% our counters for theorems, comments, remarks, etc
\newtheorem{thrm}[theorem]{Theorem}
\newtheorem{cor}[theorem]{Corollary}
\newtheorem{rmk}[theorem]{Remark}
\newtheorem{defn}[theorem]{Definition}
\newtheorem{defns}[theorem]{Definitions}
\newtheorem{ctrex}[theorem]{Counterexample}
\newtheorem{lemma}[theorem]{Lemma}
\newtheorem{note}[theorem]{Note}
\newtheorem{ex}[theorem]{Example}
\newtheorem{conj}[theorem]{Conjecture}
\newtheorem{process}[theorem]{Process}
\newtheorem{openprob}[theorem]{Open Problem}

% make a new symbol, since somehow this form makes it a solid square
%\renewcommand{\qedsymbol}{\Box}.
%\newsavebox{\proofbox} \sbox{\proofbox}{\rule{17pt}{17pt}}

% Front matter (title page, dedication, etc.).
% Chris Lomont PhD Thesis 2003
%
%  front.tex     June 3, 2002     Mark Senn
%
%  This is the ``front matter'' for a simple, example thesis.
%
%  Regarding ``References'' below:
%      KEY    MEANING
%      PU     ``A Manual for the Preparation of Graduate Theses'',
%             The Graduate School, Purdue University, 1996.
%      TCMOS  The Chicago Manual of Style, Edition 14.
%      WNNCD  Webster's Ninth New Collegiate Dictionary.
%

  % Title page is required.
  % The title page is constructed using the information specified
  % with \title, \author, \degree\ and \majorprof earlier.
\maketitle

  % Dedication page is optional.
  % A name and often a message in tribute to a person or cause.
  % References: PU 15, WNNCD 332.
\begin{dedication}
To dad, whose interest in science inspired mine,

and mom, whose support and encouragement made this possible
\end{dedication}

  % Acknowledgements page is optional but most theses include
  % a brief statement of appreciation or recognition of special
  % assistance.
  % Reference: PU 16.
\begin{acknowledgments}
I would like to express my sincere thanks to my advisor, Professor
T.T. Moh, for the guidance, encouragement, and support during my
coursework and research. I also wish to thank Professors William
Heinzer, Brad Lucier, and Kenji Matsuki, for serving on my
committee and for discussing all the countless questions and
topics I brought to them. I especially thank Kenji Matsuki for
helping me to learn Algebraic Geometry from Hartshorne's book,
which made this thesis possible.

I thank my parents, Mary Lou and Kent, for encouraging me into
science, and for all the creative things they stimulated me with
growing up. I also thank them for allowing all the crazy ideas I
pursued, all the experiments I performed on household items, and
all the fun things my siblings and I got to do as children.

I thank my many friends and peers at Purdue. Among them I thank
Roy Osawa for early guidance, Mark Rogers for rooming with me 4
years and all the work we put in learning aspects of ring theory,
Chris Mitchell for many topics of discussion, and Charles Crosby
for watching countless movies with me. I thank Jimmy Chen, who has
been a student of Professor Moh along with me, and who has
attended and talked at our many seminars on aspects of coding
theory, cryptography, and other areas of math.

I thank my officemate and friend of 7 years, Majid Hosseini, who
went through the entire process with me, inspired me, listened to
me, pushed me, and spent many hours talking math, hiking, and
enjoying life with me.

These last few years I especially thank Melissa Wilson, wonderful
girlfriend, fiancee, and soon to be wife, who has given me reason
to work harder, try to achieve more, and who causes my life to be
immensely richer.

Finally, I would like to thank all the video game developers
worldwide, who have given me plenty of distraction my many years
at Purdue. In particular, the makers of Heroes of Might and Magic
I, II, and III, the makers of Doom and Quake, and the makers of
various online chess games. Without them I would have graduated
much sooner and missed meeting some people important to me.

\end{acknowledgments}

  % The preface is optional.
  % References: PU 16, TCMOS 1.49, WNNCD 927.
%\begin{preface}
%\end{preface}

  % The Table of Contents is required.
  % The Table of Contents will be automatically created for you
  % using information you supply in
  %     \chapter
  %     \section
  %     \subsection
  %     \subsubsection
  % commands.
  % Reference: PU 16.
\tableofcontents

  % If your thesis has tables, a list of tables is required.
  % The List of Tables will be automatically created for you using
  % information you supply in
  %     \begin{table} ... \end{table}
  % environments.
  % Reference: PU 16.
\listoftables

  % If your thesis has figures, a list of figures is required.
  % The List of Figures will be automatically created for you using
  % information you supply in
  %     \begin{figure} ... \end{figure}
  % environments.
  % Reference: PU 16.
\listoffigures

  % List of Symbols is optional.
  % Reference: PU 17.
%\begin{symbols}
%  $m$& mass\cr
%  $v$& velocity\cr
%\end{symbols}

  % List of Abbreviations is optional.
  % Reference: PU 17.
%\begin{abbreviations}
%  abbr& abbreviation\cr
%  bcf& billion cubic feet\cr
%  BMOC& big man on campus\cr
%\end{abbreviations}

  % Nomenclature is optional.
  % Reference: PU 17.
%\begin{nomenclature}
%  Alanine& 2- acid\cr
%  Valine& 2-Amino-3-acid\cr
%\end{nomenclature}

  % Glossary is optional
  % Reference: PU 17.
%\begin{glossary}
%  chick& female, usually young\cr
%  dude& male, usually young\cr
%\end{glossary}

  % Abstract is required.
  % Note that the information for the first paragraph of the output
  % doesn't need to be input here...it is put in automatically from
  % information you supplied earlier using \title, \author, \degree,
  % and \majorprof.
  % Reference: PU 17.
\begin{abstract}
Error correcting codes are defined and important parameters for a
code are explained. Parameters of new codes constructed on
algebraic surfaces are studied. In particular, codes resulting
from blowing up points in $\proj^2$ are briefly studied, then
codes resulting from ruled surfaces are covered. Codes resulting
from ruled surfaces over curves of genus 0 are completely
analyzed, and some codes are discovered that are better than
direct product Reed Solomon codes of similar length. Ruled
surfaces over genus 1 curves are also studied, but not all classes
are completely analyzed. However, in this case a family of codes
are found that are comparable in performance to the direct product
code of a Reed Solomon code and a Goppa code. Some further work is
done on surfaces from higher genus curves, but there remains much
work to be done in this direction to understand fully the
resulting codes. Codes resulting from blowing points on surfaces
are also studied, obtaining necessary parameters for constructing
infinite families of such codes.

Also included is a paper giving explicit formulas for curves with
more \field{q}-rational points than were previously known for
certain combinations of field size and genus. Some upper bounds
are now known to be optimal from these examples.
\end{abstract}

\blankpage
 % todo - include

% Introduction (first chapter).
% math PhD thesis Chris Lomont 2003
%
%  intro.tex     June 3, 2002     Mark Senn
%
%  This is the introduction chapter for a simple, example thesis.
%

\chapter{INTRODUCTION} % todo - can I boldface the table of contents chapters?

\qtn{Richard W. Hamming}{Mathematics is an interesting
intellectual sport but it should not be allowed to stand in the
way of obtaining sensible information about physical processes.}

Redundant information is almost a necessity in any area of
communication. For example, the text ``CN U RD THS? THN U R ERRR
CRRCTNG" is still readable to most people, who process these
words, mentally adding missing letters as necessary. ISBN numbers
on books are engineered to detect if any digit is incorrect, and
to detect any transposition of digits. UPC codes are designed to
detect errors. Many areas of natural communication use error
correction, and certainly technological communication would be
impossible without it.

Transmitting information is one of the cornerstones of modern
technology. Due to inherent noise, this information must be
protected against corruption, and the most obvious way is to send
multiple copies of the data, then choosing the most common
outcome, a method called the ``repetition code".

The disadvantage of repetition codes is that the overhead of
sending multiple copies of the information becomes unacceptably
high. While working at Bell Labs in the USA in 1948, Dr. Claude
Shannon published ``A Mathematical Theory of Communication,"
starting information theory by showing that it was possible to
encode information such that the overhead was minimal
\cite{Shannon48}. Unfortunately his proof was not constructive,
leaving room for later researchers to seek these promised codes.

Two years later, Richard Hamming, also at Bell Labs, constructed
the simple codes bearing his name \cite{Hamming50} that take four
data bits, add three check bits, and allow the correction of any
single bit error. The repetition code above would require twelve
bits to do this, so Hamming codes began the realization of more
efficient coding methods.

Besides Shannon and Hamming, many of the pioneers in information
theory worked at Bell Labs too: Berlekamp, Gilbert, Lloyd,
MacWilliams, and Sloane.

As these events unfolded, John Leech at Cambridge created similar
codes while working on group theory. His codes were based on the
remarkable 24-dimensional Leech Lattice, and were related to
sphere packings (which have become essential to coding theory).
This lattice was also important to the classification of finite
simple groups.

The most widely used class of error correcting codes, Reed-Solomon
codes, were introduced by Irving S. Reed and Gus Solomon in a 1960
paper entitled, ``Polynomial Codes Over Certain Finite Fields,"
while they were on staff at the Massachusetts Institute of
Technology's Lincoln Laboratory \cite{ReedSolomon60}. These codes
have good properties for many small fields used in practice, and
have efficient decoding algorithms, making them indispensable to
engineers for the past several decades.

The utility of error correcting codes for information transmission
was immediately apparent, and was used by NASA in all space
programs. For example:

In 1965 the NASA Mariner probe took photos of Mars at 200x200
resolution in 64 grey levels (6 bits), transmitted 8 bits per
second, and required 8 hours to transmit a picture
\cite{Posner68}. The Mariner probe in 1969 used a Hadamard code
with a rate of 6/32, and was able to correct up to 7 errors in a
received word of 32 symbols \cite[2.1, 4.1]{vanLint}. From 1969 to
1973 NASA used a binary (32,64,16) Reed-Muller code which could
correct 7 out of 32 bits, detect 8 errors, using 6 data bits and
26 check bits. Transmission was increased to 16000 bits per
second. In January 1972 Mariner 9, the first spacecraft to orbit
another planet, used this code on 600x600 pixel pictures, taking
in 100,000 bits per second at the camera, so had to store pictures
for later transmission \cite[Ex 4.2.2]{Roman}\cite{Posner68}. In
1976, Viking landed on Mars, and took color pictures \cite{Roman}.
From 1979-1981, Voyager spacecrafts took color Jupiter and Saturn
pictures, using a 4096 symbol alphabet (4 bits each of red, green,
and blue: $2^{12}=4096$) to send color, using a binary (24,4096,8)
Golay \cite{Golay49} code, which is 3 error correcting and 4 error
detecting, with rate R=12/24 \cite[Ex 4.2.3]{Roman}

Modern consumer devices are riddled with error correcting schemes.
Modems on computers use codes to fight phone line noise. Compact
Discs (CDs) use a Cross Interleaved Reed Solomon Code (CIRC) to
protect against scratches, dirt, and cracks, and will correct up
to about 4000 consecutive errors (about 2.5 mm of track). Audio
systems can overcome even more damage by interpolating the signal
\cite{Peek85}. Computer hard drives use a Reed Solomon code to fix
errors on platters. DVDs, satellites, fax machines, and
telecommunications equipment all use error correcting codes. DARS
(digital audio radio services) and SDARS (satellite digital audio
radio services) rely on error correction to create CD quality
radio. For example, the new subscription digital radio band, XM,
at 2332.5 to 2345 MHZ, has 50 CD quality (64kb/s) channels and
many lower quality channels, and uses a Reed Solomon outer code
with a 1/2 convolutional inner code\cite{XM}. HDTV (High
Definition Television) has error correcting codes built into the
specifications. A Yahoo search on ``error correcting code" yielded
14900 hits; searching for ``coding theory" yielded 33,200 hits
(Mar 2003). Researchers think perhaps DNA uses error correction to
avoid fatal defects \cite{DNACode}. Internet transmission uses
error correction at many levels.

Perhaps the most interesting area where error codes are being
applied is in quantum computing. Due to the very sensitive nature
of quantum states, quantum computers have very special
requirements to maintain data integrity. Numerous researchers have
worked on error correcting codes to make quantum computation
feasible, \cite{Shor96, CRSS97, CRSS98}, and it was while studying
quantum computing that I became interested in error correcting
codes.

This thesis is concerned with constructing new classes of error
correcting codes, deducing their parameters, and finding better
codes usable in the future, since current codes are losing some of
their usefulness as data rates increase and demands become more
stringent. Theoretically, the best class of codes currently are
the Goppa codes, which are codes from linear systems on algebraic
curves. Their usefulness comes from the many \field{q} rational
points on curves of large genus and some deep properties about
modular curves. Since higher dimensional varieties would have even
more \field{q} points, I wanted to mimic for surfaces some of the
constructions on linear systems from curves, to see if I could
obtain better codes.

In particular, by examining codes coming from ruled surfaces I was
able to construct some codes that are better than the direct
product of Reed Solomon codes over a fixed field. This was done by
classifying all such codes from surfaces ruled over $\proj^1$.
Also codes on surfaces ruled over an elliptic curve are studied,
and partial results are obtained, giving another class of codes
that have good parameters.

Blowing up points on surfaces to obtain long families of codes is
briefly studied, but turned out to be a difficult path to analyze.
Necessary conditions are found to construct such families of
codes.

Finally, new curves are explicitly given that have more \field{q}
rational points on them than were previously known for certain
genus and $q$ combinations. In some cases this increases the
number of points to match known bounds, showing that the bounds
cannot be improved for those combinations.

The layout of this thesis is as follows. Chapter 2 contains basic
error correcting coding background. Chapter 3 covers some theorems
giving sufficient conditions to construct codes on higher
dimensional varieties. Chapter 4 is an initial attempt to
construct codes on families of surfaces by blowing up points.
Chapter 5 covers constructions on ruled surfaces, and in
particular, constructs 2 families of codes with explicit
parameters. This is done by classifying all codes on surfaces
ruled over $\proj^1$, which gives codes slightly better than the
product code of two Reed Solomon codes. The other family is over
certain surfaces ruled over an elliptic curve, and these codes are
also comparable to the corresponding product code. Chapter 6
contains new curves with many rational points, and compares them
to bounds on the number of rational points. Chapter 7 is the
conclusion and lists some open problems.

% end of file
 % todo - include

% Put other chapter comments here.
% Put other chapter \include's here.
\todo{quotes?}

\todo{ensure all blank pages have "This page deliberately left
blank" on them}

\chapter{CODING THEORY BACKGROUND}
% coding theory background - Chris Lomont PhD 2003
\section{Error Correcting Codes}
An error correcting code is a method of adding redundancy to data,
so that if the resulting redundant data gets corrupted, the
original data can be reconstructed from the corrupted data. A
linear error correcting code (LECC) is a subspace C (the
\emph{codewords}) of a vector space V over some finite field
\field{q}. Let $n=\dim V$ (the \emph{length} of the code), $k=\dim
C$ (the \emph{dimension} of the code), and then denoting
$c=(c_1,c_2,\dots,c_n)\in C$, let $d=\min\{\# c_i\neq 0 \,\, |
\,\, c\neq 0 \in C\}$ (the \emph{distance} of the code). $C$ is
called an $[n,k,d]$ code. The vector space $\field{q}^k$ (elements
are called \emph{messages}) can be embedded in V with image C.
This embedding adds the redundancy needed to correct errors. Given
a vector $m \in C$ (a codeword), and arbitrarily changing at most
$t=\lfloor \frac{d-1}{2}\rfloor$ entries in $m$ to get an element
$r\in V$ (the received word), then $m$ is the unique codeword in
$C$ that differs from $r$ in at most $t$ places. Given $r\in V$,
finding such $m\in C$ is called decoding. Creating efficient
decoding algorithms is usually separate from constructing good
codes.
\begin{note}
\emph{For this thesis, the word ``code" will denote a linear error
correcting code, unless otherwise stated.}
\end{note}

\section{Parameters of a Code}
Important parameters of an error correcting code are the rate
$R=\frac{k}{n}$ and the relative error correcting capability
$\delta=\frac{d}{n}$. Both values are in $[0,1]$, and both are
desired to be as large as possible. Of course as one value
increases, there are bounds forcing the other value to decrease.
It is a hard open problem to understand completely this
relationship for general codes.

\subsection{Shannon's Noisy Coding Theorem} Shannon's Noisy Coding
Theorem \cite{Shannon48} says (roughly) that given any rate $R$
less than the maximal rate a channel can support (which we leave
undefined for this thesis, however see \cite{Roman}), and any
decoding failure probability $p>0$ desired, there exists a code
with rate $R$ and probability of failure $<p$ (which depends on
$\delta$ and the channel noise), if \emph{the length of the code
is allowed to grow arbitrarily}. Explicit construction of such
codes is unknown; since such codes must be very long, long
families of codes are often studied. Finding codes promised by
Shannon's Theorem is a central problem in coding theory. As a
result of this theorem, people are led to look at families of
codes, whose lengths tend to infinity.

Thus the definition:
\begin{definition}\label{d:goodfamily}
A \emph{\textbf{good family}} of codes is a sequence
$C_1,C_2,\dots,C_n,\dots$ of codes over a fixed \field{q} such
that both $\limsup \delta(C_n)$ and $\limsup rate(C_n)$ are
bounded away from 0, and $\lim\limits_{n\rightarrow\infty}
length(C_n) = \infty$.
\end{definition}

For a fixed $\delta\in[0,1]$, finding a good family of codes with
the highest asymptotic rate $R$ is an important theoretical
question, but is in general unsolved. This highest value is often
denoted $\alpha(\delta)$ or $\alpha_q(\delta)$ (over an alphabet
with $q$ elements) (see \cite[Ch. 5]{vanLint}).

\subsection{Upper Bounds}
There are relations for possible values of $n$, $k$, and $d$. For
example, for an $[n,k,d]$ code, the Singleton Bound is \cite[Cor
5.2.2]{vanLint})
\begin{equation}\label{e:singleton}n+1\geq k+d
\end{equation}
The singleton bound leads to an upper bound of $\alpha(\delta)\leq
1-\delta$. There are many bounds improving on this, see for
example \cite{Roman} and \cite{vanLint}.

\subsection{Lower Bounds}
Lower bounds are important since they guarantee existence of a
family of codes. One of the best bounds is the Gilbert-Varshamov
(GV) bound \cite[Theorem 4.5.26]{Roman}:

\begin{theorem}[Gilbert-Varshamov Bound]

For an alphabet of $q$ symbols (e.g., a finite field \field{q}),
let $\theta = \frac{q-1}{q}$. \textbf{If} $0\leq\delta\leq\theta$
\textbf{then}
$$\alpha_q(\delta)\geq 1+\delta \log_q\delta +
(1-\delta)\log_q(1-\delta)-\delta\log_q(q-1)$$
\end{theorem}

Their method is constructive, yet not efficient for
implementation, thus leaving room to find efficient codes. However
for about 30 years, researchers doubted this bound could be
improved, and were surprised when it was surpassed by Goppa codes
(see section \ref{s:curvecodes} below).

\section{Direct Product Codes}\label{s:productcodes}

Let $C_i$ be an $[n_i,k_i,d_i]$ code with rate $R_i$ and relative
distance $\delta_i$, for $i=1,2$, over the same field \field{q}.
Since the code $C_i$ takes a vector of length $k_i$ and encodes it
into a vector of length $n_i$, it is natural to define a product
code $C=C_1\times C_2$ as follows: take as a message a $k_1\times
k_2$ matrix over \field{q}. Encode each row using the code $C_1$
to get a $n_1\times k_2$ matrix, then apply $C_2$ to encode each
row, resulting in a codeword in $C$, which is viewed as a
$n_1\times n_2$ matrix. With setup we have
\begin{theorem}
Given $[n_i,k_i,d_i]$ codes $C_i$ with rates $R_i$ and relative
distances $\delta_i$, $i=1,2$, the \textbf{direct product code}
$C=C_1\times C_2$ is an $[n_1n_2,k_1k_2,d_1d_2]$ code. In
particular it has rate $R_1R_2$ and relative distance $\delta_1
\delta_2$.
\end{theorem}
\begin{proof}

See \cite[ex 3.8.12]{vanLint}. The proof is an exercise, with a
solution in the back of the book which would take us too far
afield.
\end{proof}

\begin{defn}
For the rest of this paper direct product codes will be called
merely ``\textbf{product codes}".
\end{defn}

Since rates and relative distances are in $[0,1]$, product codes
never increase the rates or error capabilities. However, there are
often other reasons to mix codes: for example, to fight long burst
errors in CD players and space satellites, or to combat other
types of errors due to engineering constraints. For methods of
combining codes to obtain new codes see \cite[Ch 4.3]{Roman}.

\section{Codes from Curves}\label{s:curvecodes}
Probably the most commonly used codes are the Reed Solomon (RS)
codes, which, over a finite field \field{q} of $q$ elements, give
codes with parameters $[q-1,k,q-k]$, $0\leq k \leq q-1$. Note
these meet the Singleton Bound (equation \ref{e:singleton}), but
suffer from being a fixed length for a given field. In 1981 Goppa
\cite{Goppa} generalized RS codes to codes on algebraic curves,
with the RS codes being the Goppa code over the curve $\proj^1$.
Using the Riemann-Roch theorem \cite[10.5.1]{vanLint} to deduce
the parameters, a curve of genus $g$ with $n$ distinct
\field{q}-rational points gives $[n-1,m-g+1,n-m]$ codes, for
$2g-2<m<n$.

Briefly, let $X$ be a nonsingular curve of genus $g$ over
\field{q}, with distinct \field{q}-rational points
$P_0,P_1,\dots,P_n$. Let divisors $D=P_1+P_2+\dots+P_n$ and
$G=mP_0$, with $2g-2<m<n$, and let $K$ be the canonical divisor.
Let the code $C$ be the image of $\alpha:\sheaf{L}(G)\rightarrow
\field{q}^n$ given by $\alpha(f)=(f(P_1),f(P_2),\dots,f(P_n))$. If
$f\in \ker\alpha$ then $f$ has $\geq n$ zeros, but the number of
poles of $f$ is $\leq m < n$ which implies $f=0$. Thus $m<n
\Rightarrow \alpha$ is injective. Riemann Roch gives $\dim(G) -
\dim(K-G)=\deg(G)-g+1 = m-g+1$ so the dimension $k=\dim(G)$ of the
code is at least $m - g + 1 + \dim(K-G)$. When
$\dim(K)=2g-2<m=\deg(G)$, this reduces to $d\geq m-g+1$.

Suppose $\alpha(f)$ is nonzero, with $d$ nonzero entries. Then $f$
has $n-d$ zeros among the $P_i$, so has a pole of order at least
$n-d$ at $P_0$, forcing $\deg G \geq n-d$, so the bound on the
distance is $d \geq n-\deg G=n-m$. This guarantees a
$[n,m-g+1,n-m]$ code.

Bounds on the length of such codes are discussed in my paper
\cite{Lomont}, which is reproduced in chapter \ref{c:newcurves}.

\begin{note}\label{RR}
\emph{The Riemann Roch theorem will be used throughout this paper,
as stated in \cite[IV, Theorem 1.3]{Hartshorne} or
\cite[10.5.1]{vanLint}}
\end{note}

\begin{note}
\emph{Unless otherwise stated, all varieties and bundles in this
paper will be defined over a finite field \field{q}. In particular
every variety is assumed to have at least one \field{q} rational
point.}
\end{note}

Surpassing the GV bound was thought impossible for over 30 years,
until 1982, when Tsfasman, Vl\u{a}du\c{t}, and Zink \cite{TVZ}
used modular curves \cite{Moreno} to construct Goppa codes
surpassing the GV bound, giving the TVZ bound:

\begin{theorem}[TVZ Bound]\label{t:TVZbound}Fix a finite field \field{q}.
Let $\gamma = (\sqrt{q}-1)^{-1}$. \textbf{Then}
$$\delta+\alpha_q(\delta)\geq1-\gamma$$
\end{theorem}

Drinfeld and Vl\u{a}du\c{t} \cite{DV83} showed that the TVZ bound
is the best possible using Goppa codes.

\begin{figure}%[placement]
\includegraphics*{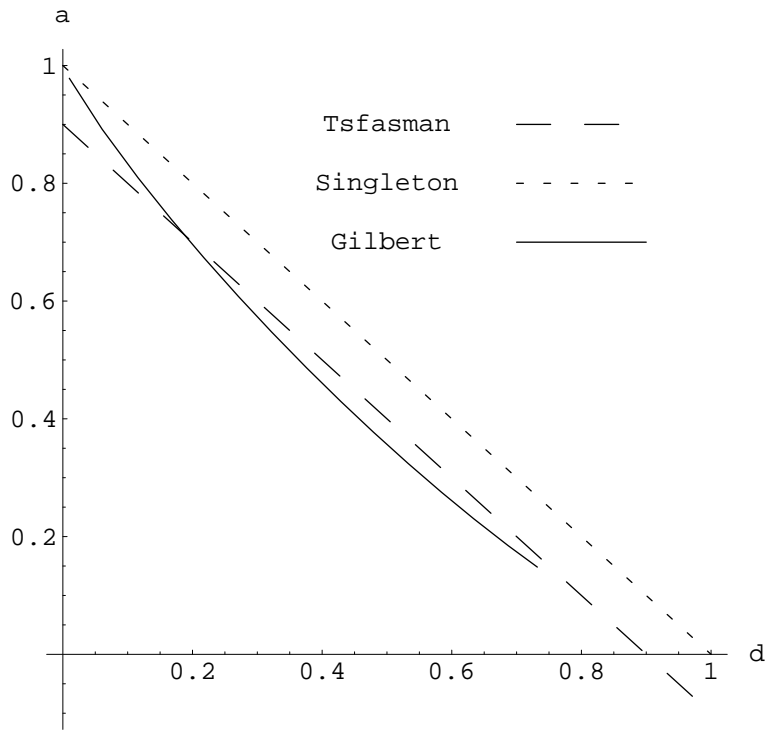}
\caption{Bounds on Parameters of Codes}\label{Fi:bounds}
\end{figure}

See Figure \ref{Fi:bounds} for a graph showing the Singleton,
Gilbert-Varshamov, and TGV bounds when $q=121$. It shows a graph
of the requested relative distance $d=\delta$ versus the
asymptotic rate $a = \alpha(\delta)$ for an infinite family of
codes. Note for some values of $d$ that the TVZ bound exceeds the
GV bound.

\section{Codes from Higher Dimensional Varieties}\label{s:highconst}
Tsfasman \cite{TV} generalized the Goppa construction to arbitrary
varieties as follows:

\begin{defn}[Code definition]\label{d:codedef}
\emph{Let $X$ be a normal projective variety over a finite field
\field{q},  and let $L$ be a line bundle on $X$ also defined over
\field{q}. Given $P_1, P_2, \dots ,P_n$ distinct
\field{q}-rational points on $X$, fix isomorphisms
$L_{P_i}\cong\field{q}$ at each stalk. Define the \textbf{code
$C(X,L)$} as the image of the germ map
$$\alpha:\Gamma(X,L)\rightarrow\bigoplus\limits_{i=1}^n
L_{P_i}\cong\field{q}^n$$ This map is evaluation of a section at
each $P_i$, and gives a vector space over \field{q}. }\end{defn}

The main problems are making sure $\alpha$ is injective and
estimating the distance of the code. To compute the distance, we
need to know: given $s\neq 0\in \Gamma(X,L)$, how many zeros does
$s$ have among the $P_i$? One tool to approach this question is
intersection theory, as in \cite{Fulton}.

For curves, these bounds are straightforward to derive using
Riemann-Roch as shown above, but this approach fails for higher
dimensional varieties.

The only work I know of studying codes from higher dimensional
varieties in some depth is the 1999 thesis of S. Hanson
\cite{Hanson}. Using theorems from this I analyzed several classes
of codes on surfaces, and classified some explicit cases. See the
last few paragraphs of the introduction for a synopsis of what
will follow.

\begin{note}
\emph{For the rest of this paper, $\alpha$ will be the germ map,
NOT the function $\alpha(\delta)$ relating the relative rate $R$
and relative distance $\delta$ for families of codes.}
\end{note}

\blankpage
% end of file
 % todo - include

\chapter{THEORY ON VARIETIES}
% theorems for inclusion in thesis
\section{Introduction}

Here are reproduced theorems that give the background for the work
that follows. Codes will be constructed by the method in
definition \ref{d:codedef}. In particular, sufficient conditions
are found guaranteeing that the germ map is injective, and bounds
are obtained on the distance of the resulting codes. Also, minor
corrections and changes from \cite{Hanson} are stated.

\section{Theorems}
This first theorem gives a bound on the distance of codes on
higher dimensional varieties.
\begin{theorem}\cite[Theorem 5.9]{Hanson}
\label{basetheorem} Suppose $X$ is a normal and projective variety
over \field{q}, $\dim X\geq 2$, and $C_1,C_2,\dots,C_\pts$ are
irreducible curves on $X$ with \field{q}-rational points
$P_1,P_2,\dots,P_n$. Assume there are $\leq N$ points on each
$C_i$. Let $L$ be a line bundle such that $L.C_i\geq0$ for all
$i$. Let
$$l=\sup\limits_{s\in\Gamma(X,L)}\#\{i : C_i \subseteq Z(s)\}$$
where $Z(s)$ is the divisor of zeros of $s$.

\textbf{Then} the code $C(X,L)$ has length $n$ and minimum
distance

$$d\geq n-lN-\sum_{i=1}^\pts L.C_i$$

\textbf{If} $L.C_i=\intpts\leq N$ for all $i$ \textbf{then}
$$d\geq n-lN-(\pts-l)\intpts$$
\end{theorem}

\begin{proof}
Let $s\in \Gamma(X,L)$. Let $D$ be its divisor of zeros.
$\alpha(s)\in \field{q}^n$ has $\#(D\cap\cup_i C_i)$
\field{q}-zero coordinates. $D\cap (\cup_i
C_i)=(\cup_{C_i\subseteq D}C_i)\cup(D\cap\cup_{C_i\nsubseteq
D}C_i)$, where the last intersection is proper. So $\alpha(s)$ has
at most $lN+\sum_{C_i\nsubseteq D}L.C_i$ zeros. $L.C_i\geq 0$, so
the last formula is bounded by $\sum_{i=1}^\pts L.C_i$, implying
$$d\geq n-lN-\sum L.C_i$$
If each curve counts the same in the intersection product ($L.C_i
= \intpts$), then we can correct for double counted zeros by
subtracting $l \intpts$ from the possible number of zeros: $$d\geq
n - lN - (\pts-l)\intpts$$
\end{proof}

\begin{cor}\cite[Cor 5.10]{Hanson}
\textbf{If} $n > lN+\sum_i L.C_i$ \textbf{then} $\alpha$ is
injective.
\end{cor}

\begin{proof}
The distance $d>0$ implies injectivity.
\end{proof}

\begin{cor}\cite[Cor 5.11]{Hanson}
\label{boundcorollary} \textbf{If} $X$ is a nonsingular surface,
$H$ is a \emph{nef} divisor on X with $H.C_i
> 0$ \textbf{then} $$l\leq\frac{L.H}{\min_i\{C_i.H\}}$$ Thus
\textbf{if} $L.H < C_i.H$ for all $i$, \textbf{then} $l=0$ and
$d\geq n-\sum_{i=1}^\pts L.C_i$
\end{cor}

\begin{note} \emph{$H$ nef (numerically effective) means $H.C\geq 0$ for all curves $C$ on
$X$.}\end{note}

\begin{proof}
Let $D$ be a member of the linear system L corresponding to
covering $l$ of the $C_i$. Then $H$ nef $\Rightarrow L.H=D.H\geq
\min\{C_i.H\}l\Rightarrow$
$$l\leq\frac{L.H}{\min\{C_i.H\}}$$
\end{proof}

The above proves

\begin{theorem}[Main Theorem]\cite[Theorem 5.1]{Hanson}
\label{maintheorem} Let $X$ be a nonsingular projective surface
over \field{q}. Let $C_1,C_2,\dots,C_\pts$ be irreducible curves
on $X$ with \field{q} rational points $P_1,P_2,\dots,P_n$. Let $L$
be a divisor on $X$ with $L.C_i\geq 0$ for each $i$. Let $H$ be a
divisor on $X$ so that H is nef and $H.C_i
> 0$. Assume $L.H < C_i.H$ for all $i$.

\textbf{Then} the code $C(X,L)$ has length $n$, minimum distance
$d\geq n-m$ where $m=\sum L.C_i=\sum \deg_L(C_i)$, and \textbf{if}
$m<n$, \textbf{then} the dimension of C is
$k=\dim_{\field{q}}\Gamma(X,L)$.
\end{theorem}

\begin{note}
\emph{Hanson \cite{Hanson} leaves out the word ``irreducible'' and
requires $H$ is ample and $L$ is nef. However, the above
conditions are strong enough to prove the theorem.}
\end{note}

\begin{note}
\emph{Bjorn Poonen \cite[Cor 3.5]{Poonen} has shown existence of
\emph{space filling curves} with the following:}
\begin{cor}
Let $X$ be a smooth, projective, geometrically integral variety of
dimension $m\geq 1$ over \field{q}, and let $E$ be a finite
extension of \field{q}. Then there exists a smooth, projective,
geometrically integral curve $Y\subseteq X$ such that $Y(E)=X(E)$.
\end{cor}
\end{note}

The methods in that paper \cite{Poonen} can perhaps be extended to
find coverings of surfaces by curves each with an equal number of
rational points, allowing the second bound on the distance from
theorem \ref{basetheorem} to be used.

\blankpage

% end of document
 % todo - include

\chapter{CODES FROM BLOWING UP POINTS ON SURFACES}
% blowing up points in P^2 to obtain long codes
\section{Long Codes}

In order to obtain good families of codes (definition
\ref{d:goodfamily}), one needs to look for long codes. Given a
fixed finite field \field{q} and a code over a variety $X$ with
$\dim X \geq 2$, one way to increase the length (the number of
\field{q} rational points) is by blowing up points. The problem
becomes finding divisors on each variety that satisfy theorem
\ref{maintheorem}.

\section{Naive Construction on Surfaces}
Here we attempt to take a code on a smooth projective surface $X$,
and by blowing up points and lifting certain divisors, create
longer codes on surfaces. First we review the intersection theory
on surfaces.
\subsection{Intersection Theory of Blow-Ups of Points on Surfaces}\label{s:intcalc} If $\pi:\tilde{X}\rightarrow X$ is a
blowup of a surface at points $P_i$, $i=1,\dots,t$, with
exceptional divisors $E_i$ above each $P_i$, then $\Pic\,
\tilde{X}\cong \Pic\, X \oplus \mathbb{Z}^t$, with each
$\mathbb{Z}$ generated by an exceptional divisor $E_i$. $\Pic\,
\tilde{X}$ has intersection calculus $C,D\in \Pic\, X\Rightarrow
\pi^*C.\pi^*D=C.D$, $(\pi^*C).E_i=0$, $E_i^2=-1$, $E_i.E_j = 0$
for $i\neq j$ \cite[Ch V, Prop 3.2]{Hartshorne}.

\subsection{An Attempt at a Family of Codes}
Now let $\dots\rightarrow X_i\rightarrow X_{i-1}\rightarrow \dots
\rightarrow X_1\rightarrow X_0$ be a sequence of smooth surfaces
defined over \field{q}, with $\pi_i:X_{i+1}\rightarrow X_i$ the
blowup of $t_i$ \field{q}-rational points on $X_i$ ($t_i>0$),
which we can specify later. Subscripts on symbols in this chapter
will associate those symbols with the surface $X_i$ or with the
code on $X_i$. $X_i$ has $n_i=\#X_i(\field{q})$ rational points.
Assume $n_i>0$ for all $i$. Let $E_i^k$ be the set of exceptional
divisors from the blowup $\pi_i$, with index $k=1,2,\dots,t_i$.
Blowing up $t_i>0$ points on surface $X_i$, we obtain $q+1$ points
over each blown up point, thus $n_{i+1}=n_i+t_i(q+1)$.

We then want to use theorem \ref{maintheorem} to construct codes
on each $X_i$, so we start by covering each with curves. Cover the
$n_0$ \field{q} rational points on $X_0$ with curves
$C_0^1,C_0^2,\dots,C_0^{s_0}$. To cover all the rational points on
each $X_i$, define curves iteratively as follows. The \field{q}
rational points of surface $X_i$ will be covered by $s_i$ curves.
Assume of the $t_i$ points blown up by $\pi_i$ that $\lambda_i^j$
of them lie on curve $C_i^j$. $\sum_j \lambda^j_i\geq t_i$
(blowing up an intersection of two curves causes the inequality).
Define, for $j=1,2,\dots,s_i$, curves $C_{i+1}^j=\pi_i^*C_i^j-\sum
E_i^\beta$, where the sum is over the $\lambda_i^j$ exceptional
divisors lying over points blown up on $C_i^j$ (that is,
$C_{i+1}^j$ is the strict transform of $C_i^j$). Add the
exceptional divisors as additional curves to cover $X_{i+1}$,
giving more curves $C_{i+1}^j=E_i^k$ for enough
$j=s_i+1,\dots,s_{i+1}$ and $k$ chosen uniquely until all $E_i^k$
are chosen. $s_{i+1}=s_i+t_i$. Thus at step $i$ the curves are
(iterated) strict transforms of an original curve $C_0^j$ on
$X_0$, or come from iterated strict transforms from some
intermediate $E_{i_0}^j$ on $X_{i_0}$, or are exceptional divisors
from the previous blow up. Note that any curve $C_i^j$ thus has
$q+1$ points if it comes from some $E_i^j$ or has
$\#C_0^j(\field{q})$ points if it comes from some $C_0^j$. To
simplify notation, $C_{i+1}^j=C_i^j$ denotes the (iterated) strict
transform when the subscripts differ.

Next give surface $X_0$ line bundles $L_0$ and $H_0$ satisfying
theorem \ref{maintheorem}. We then need line bundles $L_i$ and
$H_i$ on $X_i$ satisfying the conditions \ref{c:hc}, \ref{c:lc},
\ref{c:hl}, and \ref{c:dim} from theorem \ref{maintheorem}. From
looking at divisors on $X_i$ and $X_{i+1}$, and using the
intersection calculus in section \ref{s:intcalc}, one method to
get such line bundles is to define line bundles on $X_{i+1}$,
using $H_i$ and $L_i$ from $X_i$, by $H_{i+i}=h \pi_i^*H_i-\sum_j
E_i^j$ and $L_{i+1}=h \pi_i^* L_i-\sum_j E_i^j$, for some integer
$h>0$ to be determined below. Note that $H_i$ nef implies
$H_{i+1}$ is nef.

For the rest of the code parameters, use theorem
\ref{maintheorem}. The dimension of the code on $X_i$ is $k_i =
h^0(X_i,L_i)$. Set $m_i = \sum_{j=0}^{j=s_i}L_i.C^j_i$.
$L_i.C_i^j\geq 0$ implies $m_i\geq0$. The distance of the code on
$X_i$ is $d_i \geq n_i-m_i$ if $n_i>m_i$.

Finally, in order to obtain codes on $X_i$, from theorem
\ref{maintheorem} we require for all $i$ and $j$
\begin{eqnarray}
\label{c:hc}H_i.C_i^j>0\\
\label{c:lc}L_i.C_i^j\geq 0\\
\label{c:hl}C_i^j.H_i>L_i.H_i\\
\label{c:dim}n_i>m_i\geq 0
\end{eqnarray}

We then have using the calculus in section \ref{s:intcalc}

$$H_{i+1}.C_{i+1}^j = \left\{
\begin{array}{r@{\quad:\quad}l}
(h\pi_i^*H_i-\sum_k E^k_i).(\pi_i^*C_i^j-\sum_\beta E^\beta_k) &
C_{i+1}^j=C_i^j\\
(h\pi_i^*H_i-\sum_k E^k_i).(E_i^\beta) & C_{i+1}^j=E_i^\beta
\end{array}
\right.$$
Simplifying,
\begin{equation}\label{e:hc}
H_{i+1}.C_{i+1}^j = \left\{
\begin{array}{r@{\quad:\quad}l}
hH_i.C_i^j-\lambda_i^j & C_{i+1}^j=C_i^j\\
1 & C_{i+1}^j=E_i^\beta
\end{array}
\right.\end{equation} Similarly,
\begin{equation}\label{e:lc}
L_{i+1}.C_{i+1}^j = \left\{
\begin{array}{r@{\quad:\quad}l}
hL_i.C_i^j-\lambda_i^j & C_{i+1}^j=C_i^j\\
1 & C_{i+1}^j=E_i^\beta
\end{array}
\right.\end{equation}
\begin{eqnarray}
H_{i+1}.L_{i+1} & = &
(h\pi_i^*H_i-\sum_jE_i^j).(h\pi_i^*L_i-\sum_jE_i^j)\\ \label{e:hl}
 & = & h^2 H_i.L_i-t_i
\end{eqnarray}

\todo{check all $\pi_i$ have subscript, and check all
superscripts}

Given a code on $X_0$ so that conditions \ref{c:hc}, \ref{c:lc},
\ref{c:hl}, and \ref{c:dim} are met for $i=0$ and all $j$, induct
to find conditions meeting them for all $i$ and $j$. Assume all
four conditions are met for $i-1$ and for all $j$. Condition
\ref{c:hc} succeeds in the case \ref{e:hc} gives $H_i.C_i=1$ or
$H_i.C_i=hH_{i-1}.C_{i-1}-\lambda^j_{i-1}>0$. Since
$H_{i-1}.C_{i-1}$ is assumed positive, and can be as small as 1,
this requires $h>\lambda^{i-1}_j$, for all $i-1$ and $j$. In
particular, $E_{i-1}^j$ has $q+1$ points, forcing $h>q+1$. This
condition also suffices to ensure \ref{c:lc} holds for all $i$ and
$j$, leaving \ref{c:hl} and \ref{c:dim}.

To ensure condition \ref{c:hl}, applying \ref{e:hl} to $H_i.L_i$
repeatedly gives (for $i>0$)
\begin{equation}
H_i.L_i=h^{2i}H_0.L_0-\sum\limits_{j=0}^{i-1}(h^2)^{i-1-j}t_j\label{e:lhexp}
\end{equation}
Substituting in condition \ref{c:hl} gives
\begin{equation}
H_i.C_i>h^{2i}\left(H_0.L_0-\sum\limits_{j=0}^{i-1}\frac{t_j}{(h^2)^{j+1}}\right)
\end{equation}
Since the left hand side is positive and as small as 1, and both
sides are integers, this relation is true for all $i$ if and only
if the right hand side is $\leq 0$ for all $i>0$. This is
equivalent to
\begin{equation}\label{e:finalhl}
H_0.L_0\leq\sum\limits_{j=0}^{i-1}\frac{t_j}{(h^2)^{j+1}}
\end{equation}
Since this must hold for all $i>0$, a necessary and sufficient
condition for \ref{e:finalhl} to be satisfied, and hence necessary
and sufficient for \ref{c:hl} to be met, is that
\begin{equation} \label{e:ns1}
H_0.L_0\leq\frac{t_0}{h^2}
\end{equation}

Computing $m_i=\sum_j L_i.C_i^j$ is a bit more work. To simplify
the calculation, assume that no blown up point at any stage is an
intersection of any curves; thus $\sum\lambda_k^j=t_k$. Using
equation \ref{e:lc}, in the case the $C^j_i$ is the iterated
transform of some $C_0^j$ we have
$$L_i.C_i^j=h^iL_0.C_0-\sum\limits_{k=0}^{i-1}h^{i-1-k}\lambda_k^j$$
Summing over all $s_0$ such curves gives a contribution to $m_i$
of
$$s_0h^iL_0.C_0-\sum\limits_{k=0}^{i-1}h^{i-1-k}t_k$$
For those $C^j_i$ coming from some $E_k^j$, $0\leq k<i$ fixed, we
get a contribution of
$$t_k\left(h^{i-k-1}-\frac{h^{i-k-1}-1}{h-1}\right)$$
Summing these we obtain for $i>0$
\begin{equation}
m_i = h^is_0L_0.C_0-t_0\frac{h^i-1}{h-1} +
\sum\limits_{k=0}^{i-1}t_k\frac{h^{i-k}-2h^{i-k-1}+1}{h-1}
\end{equation}
This can be rewritten
\begin{equation} \label{e:finaldim}
m_i = \frac{h^i}{h-1}\left((h-1)s_0L_0.C_0 +
\left(1-\frac{2}{h}\right)\sum\limits_{k=0}^{i-1}\frac{t_k}
{h^k}+\frac{1}{h^i}\sum\limits_{k=0}^{i-1}
t_k-t_0\left(1-\frac{1}{h^i}\right)\right)
\end{equation}
The requirement is $0\leq m_i<n_i=n_0+(q+1)\sum_{k=0}^{i-1}t_k$ to
obtain codes. However equation \ref{e:finaldim} is difficult to
analyze, although it can be shown $m_i$ grows on order $O(h^i)$.
This means $n_i$ must grow this fast or faster, which places some
bounds on the average size of the $t_k$. Recall that $h>q+1$. For
example, if all $t_i=1$, then $n_i=n_0+(q+1)i$ is less than
\ref{e:finaldim} for large $i$, so it is impossible always to take
$t_i=1$.

This is an area needing more analysis.

\todo{see if the above can be strengthened}

\todo{analysis allowing blowups of intersection points}

\section{Conclusion}

\qtn{William Shakespeare (1564-1616)}{Though this be madness, yet
there is method in't.}

Summarizing, necessary and sufficient restrictions on $h$ to
guarantee that conditions \ref{c:hc}, \ref{c:lc}, and \ref{c:hl}
are met for all $i$ and $j$ are
\begin{equation}\label{e:finalcond}
\begin{array}{rcl}
h & >& \max_j\{q+1,\lambda_0^j\}\\
\frac{t_0}{h^2} &\geq &H_0.L_0
\end{array}
\end{equation}
Simple requirements ensuring condition \ref{c:dim} are not
completely understood, further than equation \ref{e:finaldim} and
$0\leq m_i<n_i$. If these four criteria can be met for all $i$,
then we have an infinite family of codes.

In order to make this into family of \textbf{\emph{good}} codes,
the asymptotic bounds on rates and relative distances are needed.
Unfortunately, even in the simplest cases, evaluating the
dimensions $k_i=h^0(X_i,L_i)$ seems quite hard. There might be
some way to relate $k_i$ to $k_{i-1}$ and induct, but I could
prove no such results. The analysis of the relative distance is
also difficult, but seems more likely to be understood. Also,
there are many variations on the above method, such as allowing
varying values for $h$ at each blowup or changing the form of
$H_i$ and $L_i$. It seems that blowing up points to make good
families of codes will be difficult in this generality.

\todo{check and clarify the example below}
\begin{ex}
\emph{Let $C$ be a smooth curve over \field{q} with $\gamma$
rational points. Let $X_0=C\times\proj^1$ with $n_0=(q+1)\gamma$
rational points over \field{q}. Let $C_0^j$ be copies of $\proj^1$
for $j=1,2,\dots,\gamma$, disjoint and covering $(q+1)\gamma$ of
the points. $\Pic \; X=\Pic \; C\oplus\mathbb{Z}$. Choose $L_0$
and $H_0$ ample so that $L_0.H_0>0$. Each $C_0^i$ has $q+1$ points
on it. To take the first step to $X_1$, there must be some number
$t_0\leq n_0$ of points to blow up large enough so there exists an
integer $h$ with $q+1<h\leq \sqrt{t_0}/(L_0.H_0)$. It is
sufficient to take $(q+2)(L_0.H_0)<\sqrt{t_0}$. This $h$ meets the
conditions in equation \ref{e:finalcond}. But the conditions on
$m_i$ still need checked in order to continue constructing the
family.}
\end{ex}

\chapter{CODES FROM RULED SURFACES}
% the theory for ruled surfaces
\section{Notation and Theory} Now specialize to the case of ruled
surfaces following \cite[V.2]{Hartshorne}. Fixing notation: $C$ is
a smooth curve of genus $g$ defined over \field{q}. \sheaf{E} is a
locally free sheaf of rank 2 over $C$, defined over \field{q},
corresponding to a rank 2 vector bundle $E$.

\begin{defns}{~} % this hard space is needed to cause numbering to line up
\begin{enumerate}
    \item \emph{\sheaf{E} is \textbf{decomposable} if
$\sheaf{E}\cong\sheaf{L}_1\oplus\sheaf{L}_2$ for invertible
sheaves $\sheaf{L}_i$ on $C$}
    \item \emph{\sheaf{E} is \textbf{normalized} if $H^0(\sheaf{E})\neq 0$ but
$H^0(\sheaf{E}\otimes\sheaf{L})=0$ for all invertible sheaves
\sheaf{L} on $C$ with $\deg\sheaf{L} < 0$}
    \item \emph{If \sheaf{E} is normalized define
    $e=-\deg\bigwedge^2\sheaf{E}$}
    \item \emph{$X=\proj(Symm(\sheaf{E}))$ is a ruled surface, equipped with
$\pi:X\rightarrow C$, a $\proj^1$ bundle with a section and a
relatively ample line bundle $\sheaf{O}_X(1)$. Let $C_0$ be the
corresponding divisor. \sheaf{E} normalized implies $C_0$ is an
image of a section of $\pi$}
\end{enumerate}
\end{defns}

Assume \sheaf{E} is normalized. If \sheaf{E} is decomposable then
$e\geq 0$. All such values of $e$ are possible: for example taking
$\sheaf{E}=\sheaf{O}\oplus\sheaf{O}(-e)$. If \sheaf{E} is
indecomposable then $-g\leq e\leq 2g-2$ \cite[Theorem 2.12b,ex.
2.5]{Hartshorne}. When \sheaf{E} is normalized, $e$ is an
invariant of the surface $X$.

Then we have the following facts from \cite[V Proposition
2.9]{Hartshorne}; using the notation above:
\begin{lemma}
Let $f$ be a fiber of $\pi$.

$\emph{\Pic}\, X \cong \mathbf{Z}\oplus \pi^*\emph{\Pic}\, C$
where $\mathbf{Z}$ is generated by $C_0$.

\emph{Num }$X  \cong \mathbf{Z}\oplus\mathbf{Z}$, $C_0.f = 1$,
$f^2=0$, $C_0^2=-e$.

Let $D=a C_0 + b f$ be a divisor, $p =$ char \field{q}. Define

$\kappa = \left \{ \begin{array}{ll} e & e\geq 0 \\
\frac{1}{2}\,e & e<0 \quad g<2
\\ \frac{1}{2}\,e+\frac{g-1}{p} & e<0  \quad g\geq 2\end{array}\right. $

\textbf{Then} regarding the divisor $D=aC_0+bf$: \cite[V, Theorems
2.20-2.21 and exercise 2.14]{Hartshorne}

\textbf{If} $e\geq0$ or $g\leq 1$ \textbf{then} $D$ is ample (nef)
$\Leftrightarrow a>0$ and $b>a\kappa$ (resp. $a\geq0$ and $b\geq
a\kappa$).

In positive characteristic, \textbf{if} $e<0$ and $g>1$
\textbf{then} $D$ is ample (nef) $\Rightarrow a>0$ and
$b>\frac{1}{2}a e$ (resp. $a\geq 0$ and $b\geq \frac{1}{2}a$).
$a>0$ and $b>a \kappa \Rightarrow D$ ample. \todo{nef?}

\textbf{If} $e\geq 0$ and $Y\not\equiv C_0$ is an irreducible
curve on $X$ with $Y\equiv D$ \textbf{then} $a>0$ and $b\geq a
\kappa$.

\todo{If irreducible case - Hartshorne exercise 2.14? (See notes
6-2-02 - from Hartshorne and exercises).}

\textbf{If} $E$ is the direct sum of 2 ample line bundles on $C$
\textbf{then} $C_0$ is ample \cite[Theorem 3.1.1]{Hartshorne70}
\end{lemma}

\begin{note}
\emph{\cite{Hanson} had $\frac{p}{g-1}$ instead of the correct
$\frac{g-1}{p}$ listed above.} \todo{check all notes, remarks, etc
are not italic}
\end{note}

Then the main result for ruled surfaces is \cite[Theorem
5.29]{Hanson}:

\begin{thrm}
\label{ruledtheorem} Let $C$ be a nonsingular curve of genus $g$,
\sheaf{E} a normalized vector bundle of rank 2 over $C$, and $X$
the associated ruled surface $\pi:X=\proj(S(\sheaf{E}))\rightarrow
C$, with invariant $e\geq-g$. $f$ is a fiber over a point $P_0\in
C$, and $\pts=\#C(\field{q})$. Fix integers $a\geq0$ and $b\geq0$.
If \sheaf{E} is not ample set $l=a(\lceil k \rceil - e)+b \quad
(=b$ if $e\geq 0$), else $l=b-ae$. \textbf{If} $l < \pts$ and the
bound on $d$ is positive, \textbf{then} there are $[n,k,d]$ codes
with parameters:
$$\begin{array}{l}n=(q+1)\,\pts \\
k = h^0(C,Symm^a(\sheaf{E})\otimes\sheaf{O}_C(b P_0))\\ d\geq
n-(\pts-l)a-(q+1)l
\end{array}$$
\end{thrm}

%\todo{todo - check that l=b-ae is correct, not l=b-e as Hanson
%claims}
%\todo{find ctr example}

\begin{proof}
Let $f_1,f_2,\dots,f_\pts$ be the fibers over the \field{q} points
of $C$. These disjoint lines contain all the \field{q}-rational
points of $X$, and are the curves $C_i$ in theorem
\ref{basetheorem}. Let $L\equiv aC_0+bf$. \cite[V, 2.1-4, and II,
7.11]{Hartshorne} $\Rightarrow$
$$\Gamma(X,L)\cong\Gamma(C,\pi_*L)\cong\Gamma(S^a(\sheaf{E})\otimes\sheaf{O}_C(bP_0))$$
%($a\geq1$, $b\geq a\kappa$ (TODO?)) - reason for H nef below?
Let $H=C_0+\lceil\kappa\rceil
f$. Then $H$ is nef, $H.f_i=1$, $L.f_i=a$ for all $i$, and
$$\begin{array}{rcl}
 H.L & = & aC_0^2+(b+\lceil\kappa\rceil a)C_0.f+b\lceil\kappa\rceil f^2 \\
     & = & -e a+b+a\lceil\kappa\rceil \\
     & = & a(\lceil\kappa\rceil-e)+b \\
   ( & = & b \mbox{ if } e\geq0)
   \end{array} $$
By Theorem ~\ref{maintheorem} and Corollary ~\ref{boundcorollary},
$l\leq\frac{L.H}{\min\limits_i\{C_i.H\}} =
\frac{a(\lceil\kappa\rceil-e)+b}{1}\Rightarrow$
$$l\leq a(\lceil\kappa\rceil-e)+b$$
$\sheaf{E}$ ample $\Rightarrow C_0$ nef, so let $H=C_0$, which
gives $l=H.L=b-ae$.
\end{proof}
\todo{redefine proofbox symbol in puthesis.sty to fix these solid
boxes - see Math into Latex p323}

\begin{rmk}
\emph{\cite[5.29]{Hanson} omits \emph{normalized} in the statement
above, and states that if \sheaf{E} is ample, then $l=b-e$, which
is incorrect. See sections \ref{s:cterexp1} and \ref{ex1} for
counterexamples if normalized is omitted. \cite[5.30]{Hanson} uses
a non-normalized sheaf, and computes $e$, which is then not an
invariant, and is not guaranteed to satisfy $e\geq -g$.
\cite{Hanson} does not clearly define $e$ to be the invariant, and
seems to treat it both ways. For example:}
\end{rmk}

\begin{ctrex}
\label{ex1} \emph{Over $\proj^1$, which has genus 0, the
non-normalized sheaf $\sheaf{E}=\sheaf{O}(1)\oplus\sheaf{O}(1)$
has $\deg\bigwedge^2\sheaf{E}=-2$ which violates $e\geq -g$.}
\end{ctrex}

% end of document

\section{Ruled Surfaces over $P^1$}
% the case of genus 0 curves - to be included in PhD thesis
First we classify codes on ruled surfaces over the unique genus 0
curve, $\proj^1$, using theorem \ref{ruledtheorem}.
\subsection{Codes}
All indecomposable vector bundles over $\proj^1$ are trivial
\cite{Grothendieck}. Pic $\proj^1\cong\mathbb{Z}$, generated by a
hyperplane section. Thus $\sheaf{E}\cong \sheaf{O}(t)\oplus
\sheaf{O}(u)$ for some integers $t$ and $u$, where
$\sheaf{O}=\sheaf{O}_{\proj^1}$. Since $\proj (S(\sheaf{E}))\cong
\proj(S(\sheaf{E}\otimes\sheaf{O}(n)))$ for any integer $n$
\cite[II.7, ex 7.9b]{Hartshorne}, reduce to the normalized,
decomposable case with $\sheaf{E}\cong
\sheaf{O}\oplus\sheaf{O}(-e)$ for some integer $e\geq 0$. Then
\sheaf{E} is not ample, so $l=b$ in \ref{ruledtheorem}. All
possible codes using this construction on ruled surfaces over
$\proj^1$ have the following parameters in the notation of theorem
\ref{ruledtheorem} (corresponding to a divisor $D\sim aC_0+b f$ on
$X$):
$$\begin{array}{l}
\pts=q+1 \\
n=(q+1)^2\\
k=h^0(\proj^1,S^a(\sheaf{E})\otimes\sheaf{O}(b P_0))\\
d\geq(q+1)^2-(q+1-b)a-b(q+1)
\end{array}$$
We require $b<\pts=q+1$ from theorem \ref{ruledtheorem}.

The bound on $d$ is required to be positive to ensure the germ map
$\alpha$ is injective, so we simplify the distance bound:
$$d\geq (q+1-a)(q+1-b)$$
Note this does not depend on the choice of $e$. Since $b<q+1$,
requiring $(q+1-a)(q+1-b)>0$ is equivalent to requiring $a<q+1$.
So for a given $q$ and $e$ combination, there are $(q+1)^2$ codes
under this construction, corresponding to $0\leq a,b<q+1$.

To compute the dimension $k$, note that $S^a(\sheaf{E})\otimes
\sheaf{O}(b P_0)\cong \bigoplus\limits_{j=0}^{j=a}\sheaf{O}(b-j
e)$. Over $\proj^1$, since cohomology commutes with direct sums,
Riemann-Roch gives the dimension $k=\sum(b-j e+1)$, where the sum
is over nonnegative $j$ such that $j e \leq b$ and $j\leq a$. So
we see that the dimension of the code can be increased by
increasing either a or b. However, as usual in coding theory, this
decreases the distance $d$ of the code.

To evaluate the performance of these codes, notice that increasing
$e$ decreases the dimension $k$ of the code, and leaves the
distance $d$ unchanged, so taking $e=0$ will result in the largest
dimension $k$ for fixed $a,b$. Then we have
$k=\sum\limits_{j=0}^{j=a}(b+1)=(a+1)(b+1)$. This gives

\begin{theorem}[Lomont Code \#1]\label{t:p1codes}
The construction of Theorem \ref{ruledtheorem} applied to ruled
surfaces over $\proj^1$, defined over the finite field \field{q},
results in $[n,k,d]$ codes with parameters
$$\begin{array}{l}
n = (q+1)^2\\
k = \sum (b-je+1)\\
d \geq (q+1-a)(q+1-b)
\end{array}$$
where the sum is over nonnegative $j$ so that $j e \leq b$ and $j
\leq a$, and $0\leq a,b<q+1$, $0\leq e$.

The codes with highest rates are then when $e=0$, giving codes for
integers $0\leq a,b<q+1$
\begin{equation}\label{e:bestp1codes}
\begin{array}{l}
n = (q+1)^2\\
k = (a+1)(b+1)\\
d \geq (q+1-a)(q+1-b)
\end{array}\end{equation}
In this case $X\cong\proj^1\times\proj^1$.

\end{theorem}

\begin{proof}
Shown above.
\end{proof}

As a sanity check compare this to the Singleton Bound (see
equation \ref{e:singleton}), $n+1\geq k+d$. Theorem
\ref{t:p1codes} gives
$$
\begin{array}{rcl}
(q+1)^2+1 & \geq & (a+1)(b+1)+(q+1-a)(q+1-b)\\
1 & \geq & (a+1)(b+1)-(a+b)(q+1)+ab\\
1 & \geq & 2ab-q(a+b)
\end{array}$$
and it can be shown that for $0\leq a,b\leq q$ this is always
satisfied.

\begin{note}[Decoding]\label{n:decodep1}
\emph{The Lomont code \#1 has parameters which look like a product
code. There is a way to define a Goppa code over $\proj^1$ so it
has $q+1$ points by taking the poles at a point not defined over
\field{q}. Then the decoding algorithms for the Goppa code
\cite{Pretzel} should be able to decode the Lomont code \#1 as a
product code. This needs checked.}
\end{note}

\subsection{Comparison to Product Codes}
The Reed Solomon (RS) codes over \field{q} have parameters
$[q-1,k,q-k]$ for $0\leq k \leq q-1$ \cite[Ch 8.2]{Roman}. So the
Lomont Code \#1 in theorem \ref{t:p1codes} is longer than the
product code (section \ref{s:productcodes}) from two RS codes,
which has parameters $[(q-1)^2,ab,(q-a)(q-b)]$, $0\leq a,b\leq
q-1$. For comparison, fix the field to be the commonly used field
\field{256}, and compare the best relative distance $\delta$ for a
desired rate $r$ (see Table \ref{t:code1}). The leftmost column is
the desired rate $r=0.1, 0.2,\dots,0.9$. Then all legal
combinations of $a$ and $b$ are searched to find the best
performing $\delta$ when the corresponding rate is $\geq r$. The
left block shows the optimal choices of $a$ and $b$ giving the
rate and $\delta$ combination for the Reed Solomon code, and the
right block shows the same information for the Lomont Code \#1.
For example, looking for a code with rate at least 0.8, the
highest $\delta$ product RS code is when $a=228$, $b=229$, and has
relative distance $\delta=0.1163$. The corresponding best Lomont
code is when $a=b=229$ and has $\delta=0.1187$, so is a slightly
better code. For this example the product code could correct 377
errors, but the Lomont code could correct 391 errors, and thus is
slightly better at handling burst errors. So besides being longer
than the RS product code, the Lomont code for this example has
better relative distance. Note that the Lomont code is better for
6 of the 9 values tested.

\begin{table}
\caption{Comparison of the Reed Solomon code to the Lomont code
\#1}\label{t:code1}
\begin{tabular}{|c||c|c|c|c||c|c|c|c|} \hline
& \multicolumn{3}{r}{Reed Solomon} & & \multicolumn{3}{r}{Lomont
Code \#1} & \\ \hline
rate & a & b & rate & $\delta$ & a & b & rate & $\delta$ \\
\hline \hline
0.1&81& 81& 0.1009& 0.470973& 80& 81& 0.100562& 0.47165 \\
0.2&114& 115& 0.201615&  0.307912& 114& 114& 0.20023& 0.309603\\
0.3&140& 140& 0.301423& 0.206936& 140& 140& 0.301004& 0.207255\\
0.4&161& 162& 0.401107& 0.137332& 162& 162& 0.402262& 0.136641\\
0.5&180& 181& 0.501038& 0.0876586& 181& 181& 0.501506& 0.0874502\\
0.6&198& 198& 0.602907& 0.0517339& 198& 199& 0.602583& 0.05181\\
0.7&213& 214& 0.700992& 0.0277739& 214& 215& 0.703114& 0.0273433\\
0.8&228& 229& 0.802953& 0.0116263& 229& 229& 0.800921& 0.01187\\
0.9&242& 242& 0.900638& 0.00301423& 243& 243& 0.901391&
0.00296749\\ \hline
\end{tabular}
\end{table}

\todo{include mathematica code in appendix?}

%\todo{Now suppose we want to evaluate this code to see if it
%performs well for large $q$ (giving long codes). Fix the rate at
%$0<r<1$.}

%\todo{check $\sheaf{O}(m)\oplus\sheaf{O}(n)$ case directly -
%perhaps $d$ or $k$ better.}

%\todo{show distance bound is exact in genus 0 case.}

\subsubsection{A Counterexample}\label{s:cterexp1}
If we do not require \sheaf{E} to be normalized we get impossible
codes, showing that normalized is necessary in theorem
\ref{ruledtheorem}. Let $\sheaf{F}=\sheaf{O}(t)\oplus\sheaf{O}(u)$
for an integers $t\geq u$. Then
$$S^a(\sheaf{F})\otimes\sheaf{O}(bP_0)\cong\bigoplus\limits_{j=0}^{j=a}\sheaf{O}((a-j)t+ju+b)$$
From Riemann-Roch the dimension $k=\sum at+j(u-t)+b+1$, where the
sum is over $0\leq j\leq a$ and $at+j(u-t)+b\geq 0$. Letting
\sheaf{E} be the normalized vector bundle
$\sheaf{F}\otimes\sheaf{O}(-t)\cong\sheaf{O}\oplus\sheaf{O}(u-t)$,
we have that $e=-(u-t)=t-u\geq 0$. The length of the code is
$n=(q+1)^2$, and these are the Lomont codes in theorem
\ref{t:p1codes} above.

\begin{ctrex}\label{s:cterex}
\emph{However, taking the unnormalized \sheaf{F} with $t>0$ and
$u=-t$, since all terms in the dimension sum are nonnegative, we
have that $k\geq at$. $e=0$, and \sheaf{F} is not ample, so $l=b$.
The distance is again bounded by $d\geq (q+1)^2-(q+1-b)a-b(q+1)$,
which is independent of $t$. Taking $a=1$, $b=0$, then $d\geq
(q+1)q$, so $\alpha$ is injective and we have a code. But letting
$t$ increase without bound increases $k$ without bound, a
contradiction since the unbounded $k$-dimensional vector space
cannot be a subspace of the fixed $n$-dimensional one.}
\end{ctrex}

% end of document

\section{Decomposable Bundles over Positive Genus Curves}
% the genus > 0 decomposable case
\label{s:decomposable} Let $C$ be a curve of genus $g\geq 1$ over
\field{q}, and let $\sheaf{E}=\sheaf{O}\oplus\sheaf{O}(-e)$ for
some nonnegative integer $e$, giving the resulting surface $X$ an
invariant of $e\geq 0$.

A bound on the number \pts \, of \field{q}-rational points on $C$
is $|\pts -(q+1)|\leq g\lfloor 2\sqrt{q}\rfloor$ \cite{Serre}, but
can be improved in many cases, some of which are mentioned in
chapter \ref{c:newcurves} which deals with explicit curves
reaching known bounds. Following the reasoning from the genus 0
case, we get that
$$\begin{array}{l}
n=\pts(q+1) \\
d\geq n-(\pts-b)a-(q+1)b
\end{array}$$
The dimension $k$ is estimated using
$S^a(\sheaf{E})\otimes\sheaf{O}(bP_0)\cong\bigoplus\limits_{j=0}^{j=a}\sheaf{O}(-j
e)\otimes\sheaf{O}(bP_0)$. Using Riemann-Roch, $l(D)=\deg D+1-g +
l(K-D)$ gives that $k$ is in the form $\sum(b-je)-g+1+\zeta$,
where the sum is over certain nonnegative terms as usual, and
$\zeta$ is also nonnegative from the $l(K-D)$ part of
Riemann-Roch. Since the distance does not rely on $e$, we can
maximize the dimension $k$ by making $e=0$. This proves

\begin{theorem}
Given a genus $g$ smooth projective curve $C$ over \field{q}, and
an integer $e\geq0$. Let $\pts=\#C(\field{q})$,
$\sheaf{E}=\sheaf{O}_C\otimes\sheaf{O}_C(-e)$, and $X$ be the
corresponding ruled surface. For integers $0\leq a$, $0\leq b <
\pts$ there are codes (assuming the bound on $d$ is positive)
$$\begin{array}{l}
n=\pts(q+1) \\
k= \sum(b-je)-g+1+\zeta_j\\
d\geq(q+1-a)(\pts-b)
\end{array}$$
where the sum is over those $j$ so that $0\leq j\leq a$ and
$je\leq b$. $\zeta_j=l(K-D_j)$, where $K$ is the canonical divisor
on $C$, and $D_j=\sheaf{O}(b-je)$. Again it is clear that $e=0$
gives the largest dimension, then $X\cong C\times\proj^1$, and the
best such codes are product codes.
\end{theorem}

\todo{do these meet TGV in extremes?}

\todo{other line bundles over these curves?}

Without specific curves, it is hard to go much further than this,
since one needs information about the canonical divisor $K$ and
the dimensions $\zeta_j$. An interesting case would be to study
the curves of Garcia and Stichtenoth (see section
\ref{s:GSCurves}), since they have many rational points and would
result in very long codes.

% end of file

\section{Ruled Surfaces over Elliptic Curves}
% the indecomposable over elliptic curve case

\newcommand{\edeg}{\ensuremath{\Lambda}}               % degree of indecomposable sheaf
\newcommand{\rank}{\ensuremath{\mbox{rank}\,\,}}       % rank of a sheaf

\subsection{Vector Bundles over Elliptic
Curves}\label{s:vectclass} Here recall some facts from the
classification of indecomposable vector bundles over elliptic
curves. The classification over algebraically closed fields of any
characteristic was done by Atiyah in \cite{Atiyah57}, and the
extension to perfect fields was done in the thesis of Agnes
Williams under G. Faltings,  and was stated in a paper by Arason,
Elman, and Jacob \cite{AEJ}. Thus the following also is true for
the case we need, namely the finite field case $K=\field{q}$.

The facts are for an arbitrary perfect field $K$ and elliptic
curve $C$ defined over $K$:

1. To each $K$-rational point $P\in C$ there is constructed a
vector bundle $E_{r,d}(P)$ of rank $r$ and degree $d$ on $C$.

2. Each $E_{r,d}(P)$ is shown to be absolutely indecomposable.

3. $E_{r,d}(P)\cong E_{r,d}(Q)\Rightarrow P=Q$

4. For an absolutely indecomposable vector bundle $M$ of rank $r$
and degree $d$ on $C$ there is found a rational point $P\in C$
such that $M\cong E_{r,d}(P)$.

5. There is an absolutely indecomposable vector bundle $F_r$ of
rank $r$ and degree 0 on $C$, unique up to isomorphism, such that
$F_r$ has non-trivial global sections. Moreover, there is an exact
sequence
$$0\rightarrow \sheaf{O}_C\rightarrow F_r \rightarrow
F_{r-1}\rightarrow 0$$ If $M$ is an absolutely indecomposable
vector bundle of rank $r$ and degree $d$, then there is a line
bundle $L$ of degree 0 on $C$, unique up to isomorphism, such that
$M\cong L\otimes F_r$. $M$ contains $L$ as a subbundle.

6. $\dim \Gamma(F_r\otimes F_s)=\min\{r,s\}$. Given a line bundle
$L$, $\dim \Gamma(L\otimes F_r\otimes F_s)=0$ unless $L\geq 1$
\cite[III, Lemma 17]{Atiyah57}.

7. $F_r\otimes F_s \cong \sum\limits_{j=1}^{\min(r,s)}F_{r_j}$,
$\sum\limits_j r_j=rs$ \cite[III, Lemma 18]{Atiyah57}.

\subsection{Codes}

Let $C$ be an elliptic curve over \field{q}. Let \sheaf{E} be a
rank 2 normalized vector bundle over $C$ defined over \field{q}.
The case \sheaf{E} decomposable is covered by section
\ref{s:decomposable}. If \sheaf{E} is indecomposable, then
\cite[V, Theorem 2.15]{Hartshorne} gives that $\deg\sheaf{E}$ is 0
or 1.

Using the decomposable case from section \ref{s:decomposable}
would not be difficult to do for an elliptic curve, and would
result in a product code from results in that section, so we do
not do it here.

\subsubsection{The Degree 0 Case}

We associate the (torsion-free) coherent sheaf $\sheaf{F}_r$ to
the vector bundle $F_r$. \cite{Serre55}. To compute the dimension
of the codes in this case, we need to understand the structure of
$S^n(\sheaf{F}_r)$, enough of which is given by following theorem
for our purposes.

\begin{thrm}
\label{decomptheorem} Denote by $\sheaf{F}_r$ the unique degree 0
rank r indecomposable vector bundle with a global section over the
elliptic curve $C$, both defined over the perfect field $K$.
\textbf{Then} $S^n(\sheaf{F}_r)\cong\oplus\sheaf{F}_{r_i}$, for
some $\sheaf{F}_{r_i}$, with $\sum_i r_i=\binom{n+r-1}{r-1}$.
\end{thrm}

\begin{proof}
Write \sheaf{F} for $\sheaf{F}_r$. Let
$S^n(\sheaf{F})\cong\oplus\sheaf{E}_i$, where the $\sheaf{E}_i$
are indecomposable. Let \sheaf{L} be a degree 0 line bundle with
no global section (e.g., corresponding to a divisor $P-Q$ for
$P\neq Q$). In the exact sequence
\begin{equation}\label{xseq}
0\rightarrow \sheaf{I} \rightarrow \sheaf{F}^{\otimes
n}\rightarrow S^n(\sheaf{F})\rightarrow 0\end{equation}  $\deg
\sheaf{I}=0$ since degree is additive, and the other two terms
have degree 0 by theorems \ref{tensordeg} and \ref{symmetricdeg}.
After tensoring with \sheaf{L}, this gives the cohomology sequence
\begin{equation}
\begin{array}{ccccccccc}
0&\rightarrow&H^0(\sheaf{I}\otimes\sheaf{L})&\rightarrow&
H^0(\sheaf{F}^{\otimes n}\otimes\sheaf{L})&\rightarrow&
H^0(S^n(\sheaf{F})\otimes\sheaf{L})&\rightarrow&\\
&&H^1(\sheaf{I}\otimes\sheaf{L})&\rightarrow&
H^1(\sheaf{F}^{\otimes n}\otimes\sheaf{L})&\rightarrow&
H^1(S^n(\sheaf{F})\otimes\sheaf{L})&\rightarrow&0
\end{array}\end{equation}
The higher terms all vanish  by Grothendieck Vanishing \cite[III,
Theorem 2.7]{Hartshorne}.

Then using properties 6 and 7 in section \ref{s:vectclass}, and
that for any vector bundle \sheaf{E} Riemann-Roch gives
$h^0(\sheaf{E})-h^1(\sheaf{E})=\deg\sheaf{E}$, we get that both
terms involving $\sheaf{F}^{\otimes n}\otimes\sheaf{L}$ vanish.
Thus $H^0(\sheaf{I}\otimes\sheaf{L})=0$, so using Riemann-Roch
again gives $H^1(\sheaf{I}\otimes\sheaf{L})=0$, which implies
$H^0(S^n(\sheaf{F})\otimes\sheaf{L})=0$. So $\deg\sheaf{E}_i\leq
0$ for every $i$. Since degree is additive over direct sums
\cite[II, ex 6.12(3)]{Hartshorne} it must be that
$\deg\sheaf{E}_i=0$ for every $i$.

From Property 5, section \ref{s:vectclass}, $\sheaf{E}_i\cong
\sheaf{F}_{r_i}\otimes\sheaf{L}_i$ for some degree 0 invertible
sheaves $\sheaf{L}_i$. Suppose some
$\sheaf{L}_i\ncong\sheaf{O}_C$. Then taking the sequence
\ref{xseq} and tensoring with $\sheaf{L}_i^{-1}$, we get the
cohomology sequence
\newcommand{\sli}{\sheaf{L}_i^{-1}}       % the inverse line sheaf
$$\begin{array}{ccccccccc}
0&\rightarrow& H^0(\sheaf{I}\otimes\sli) &\rightarrow&
H^0(\sheaf{F}^{\otimes n}\otimes\sli)&\rightarrow&
H^0(S^n(\sheaf{F})\otimes\sli)&\rightarrow&\\
&&H^1(\sheaf{I}\otimes\sli)&\rightarrow& H^1(\sheaf{F}^{\otimes
n}\otimes\sli)&\rightarrow&
H^1(S^n(\sheaf{F})\otimes\sli)&\rightarrow& 0
\end{array}$$
Again by the reasoning above, the terms involving
$\sheaf{F}^{\otimes n}$ vanish, causing the terms involving
\sheaf{I} to vanish, so again $H^0(S^n(\sheaf{F})\otimes\sli)=0$.
But $S^n(\sheaf{F})\otimes\sli$ has $\sheaf{F}_{r_i}$ as a direct
summand, and thus has a global section, a contradiction. So all
$\sheaf{L}_i\cong\sheaf{O_C}$, giving that
$S^n(\sheaf{F}_r)\cong\oplus_{i=1}^m\sheaf{F}_{r_i}$.

The rank of the left hand side is given in the appendix, theorem
\ref{symmetricrank}, so the right hand side, being additive, gives
$\sum_i r_i=\binom{n+r-1}{r-1}$.
\end{proof}

\begin{note}
\emph{In special cases we can find precisely which
$\sheaf{F}_{r_i}$ occur. For example, using \cite{Atiyah57} and
the techniques above, $S^2(\sheaf{F}_2)$ can be found for any
characteristic by examining possible ranks $r_i$, and counting
global sections.}
\end{note}

%\todo{the above works on any curve and sheaf that is extensions of
%\sheaf{O}?}

%\todo{remove command sli?}

Now we can compute the parameters of codes arising from
$\sheaf{F}_2$. Using theorem \ref{ruledtheorem} we obtain

\begin{thrm}[Lomont Code \#2]\label{t:genus1degree0Codes}
Let $C$ be an elliptic curve with $\pts$ \field{q} rational
points, $\pts>0$. Let $a,b$ be integers with $0<b<\pts$, $0\leq a
<q+1$. \textbf{Then} there are $[n,k,d]$ codes with
$$\begin{array}{l}n=(q+1)\,\pts \\
k = (a+1)b\\
d\geq (q+1-a)(\pts-b)
\end{array}$$
\end{thrm}

\begin{rmk}
\emph{For $b=0$ there is still a code, but the exact dimension is
slightly more complex to compute, and appears to be
uninteresting.}
\end{rmk}

\begin{proof}
We use theorem \ref{ruledtheorem}, in the case \sheaf{E} is the
unique indecomposable rank 2 degree 0 vector bundle over $C$ with
a global section. The value for $n$ follows immediately. \sheaf{E}
is not ample \cite[Theorem 1.3]{Hartshorne71}, giving invariant
$e=0$, so $l=b$ in \ref{ruledtheorem}. The dimension arises from
using theorem \ref{decomptheorem}, Riemann Roch \ref{RR}, and
theorem \ref{tensordeg} in \ref{ruledtheorem} to obtain:
$$
\begin{array}{rcl}
k & = & h^0(C,S^a(\sheaf{E})\otimes\sheaf{O}_C(b P_0))\\
& = & h^0(C,\oplus_i(\sheaf{F}_{r_i}\otimes\sheaf{O}_C(b P_0)))\\
& = & \sum_i h^0(C,\sheaf{F}_{r_i}\otimes\sheaf{O}_C(b P_0)))\\
& = & \sum_i r_i b\\
& = & \binom{a+2-1}{2-1} b\\
& = & (a+1)b
\end{array}
$$
Note that the bound on $d$ positive if and only if $a<
\frac{n-(q+1)b}{\pts-b}=q+1$, so by theorem \ref{ruledtheorem} we
have a code.
\end{proof}

\begin{note}[Decoding]\label{n:decodeg1}
\emph{Similar to note \ref{n:decodep1}, the Lomont code \#2 has
parameters which look like a product code. The methods of decoding
mentioned there should also be able to decode the Lomont code \#2.
This too needs checked.}
\end{note}

\begin{thrm}[Lomont code \#2 Rate]
Use the notation of theorem \ref{t:genus1degree0Codes}. For a
fixed relative distance $\delta=\frac{d}{n}$, set
$b_0=\pts\left(1-\sqrt{\frac{(q+1)\delta}{(q+2)}}\;\right)$, and
then $a_0=\frac{(q+1)(\pts\delta-\pts+b_0)}{b_0-\pts}$. This code
has the highest rate $R=\frac{k}{n}$ when the integers $(a,b)$ are
one of the integer lattice points points $(\lfloor a_0
\rfloor,\lfloor b_0 \rfloor)$, $(\lceil a_0 \rceil,\lfloor b_0
\rfloor)$, or $(\lfloor a_0 \rfloor,\lceil b_0 \rceil)$, depending
on which ones satisfy the requirement of the size of $\delta$.
\end{thrm}

\begin{proof}
Given a fixed value for $\delta=\frac{d}{n}$, we wish find the
values of $a$ and $b$ giving the highest rate $R=\frac{k}{n}$,
treating $a$ and $b$ as real numbers. Assume $0<b<\pts$ and
$0<a<q+1$. Then
$$
\begin{array}{rcl}
\delta & = & \frac{(q+1)\pts-(\pts-b)a-(q+1)b}{(q+1)\pts}\\
R & = & \frac{(a+1)b}{(q+1)\pts}
\end{array}
$$
Solving the first equation for $a$, and substituting into the
second,
$$
\begin{array}{rcl}
a & = & \frac{(q+1)(\pts\delta-\pts+b)}{b-\pts}\\
R & = &
\frac{b}{(q+1)\pts}\left(1+\frac{(q+1)(\pts\delta-\pts+b)}{b-\pts}\right)
\end{array}
$$
which is valid for $0<b<\pts$. The first two derivatives of $R$
with respect to $b$ are
$$
\begin{array}{rcl}
R' & = & \frac{1}{\pts(q+1)} + \frac{b^2-2\pts b+(\delta-1)\pts^2}{(b-\pts)^2\pts}\\
R'' & = & \frac{2\pts\delta}{(b-\pts)^3}
\end{array}
$$
$b<\pts\Rightarrow R'' < 0$, so roots of $R'$ give local maxima.
These roots occur when
$$
b_0 = \pts\left(1 \pm \sqrt{\frac{(q+1)\delta}{(q+2)}}\right)
$$ $b<\pts$ forces choosing the negative sign in the $\pm$,
then this gives the best possible value of $b_0$ in the theorem.
Then $a_0$ follows. Since $R$ is concave downward, the best
possible integer $b$ for any given $a$ must be $\lfloor b_0
\rfloor$ or $\lceil b_0 \rceil$. Similarly for $a$ compared to
$a_0$. Since $(a,b)$ must be integers, and $\delta$ decreases as
$a$ or $b$ increases, the combination $(\lceil a_0 \rceil,\lceil
b_0 \rceil)$ will have too small a delta value, giving the other
three combinations as possible outcomes. It is possible to
construct examples with each combination as the best choice.
\end{proof}
% see notes 9-26-02 for this above derivation - notes incorrect

%\todo{add mathematica code to appendix}

\subsection{Comparison to Product Codes}

Next we compare these codes with the product codes obtained from
Goppa codes on $C$ with Reed-Solomon codes on $\proj^1$. On
$\proj^1$ over \field{q} the RS codes have parameters
$[q-1,k_1,q-k_1]$ for $0\leq k_1 \leq q-1$. The Goppa codes on $C$
have parameters $[\pts-1,k_2,\pts-1-k_2]$ with $0<k_2<\pts-1$.
Similar to the above, for $0<k_2<\pts-1$ and $0\leq k_1\leq q-1$
we have product code parameters
$$
\begin{array}{rcl}
\delta_1 & = & \frac{(\pts-1-k_2)(q-k_1)}{(q-1)(\pts-1)}\\
R_1 & = & \frac{k_1k_2}{(q-1)(\pts-1)}
\end{array}
$$
Solving the first for $k_1$ and substituting into $R_1$ gives
$$ R_1 =
\frac{qk_2}{(q-1)(\pts-1)}-\frac{k_2\delta_1}{\pts-1-k_2}$$ The
best possible value for $k_2$ then becomes, for a fixed
$\delta_1$,
$$k_2=(\pts-1)\left(1\pm\sqrt{\frac{(q-1)\delta_1-1}{q}}\,\right)$$
We need $k_2<\pts-1$, so we take the negative sign choice.

Using the same value for $\delta$, substitute the best values to
maximize the rates for each code, and subtract the resulting
optimal rates, giving the difference between the optimal rates as
a function of $q$ and $\delta$.
\begin{equation}\label{f:err}
err(q,\delta) = R-R_1 =
-\frac{2}{q^2-1}+2\sqrt{\delta}\left(\sqrt{\frac{q}{q-1}}-\sqrt{\frac{q+2}{q+1}}
\, \right)\end{equation} Over a field size used often in practice,
$q=256$, this simplifies to
$err(256,\delta)=-0.000030518+0.0000304586\sqrt\delta$. Since
$\delta\in[0,1]$, this shows the surface code has a slightly lower
rate than the product code. Below in section
\ref{s:CompareElliptic} this is shown to be true for any size
finite field. The two rates converge to the same value as $q$ gets
larger and larger, so for large fields the Lomont code \#2 code
performs arbitrarily close to the product code from the RS code
and Goppa code.

\subsubsection{An Example}
However, since the parameters $a$, $b$, $k_1$, and $k_2$ are
restricted to integral values, sometimes the Lomont code \#2 is
slightly better than the product codes, since the rates are so
close when considered as continuous functions. For example, see
table \ref{t:code2}. Here the field is fixed at \field{256} as in
the genus 0 case, and the elliptic curve is $y^2=x^3+x+1$, which
has $\pts=255$ rational points.  The left column is the desired
rate $R=0.1,0.2,\dots,0.9$, and the rest shows the best parameters
for the product code and the Lomont code \#2. Notice in some
cases, like the $R=0.6$ case, that the Lomont code \#2 has higher
relative distance than the corresponding product code, but a lower
rate. As in the genus 0 case, since it is longer but with a
similar rate and relative distance, the Lomont code can correct
longer burst errors than the product code can. For this example,
the Lomont code \#2 can correct 1694 errors, while the product
code can correct 1659 errors.

\begin{table}
\caption{Comparison of the product code to the Lomont code
\#2}\label{t:code2}
\begin{tabular}{|c||c|c|c|c||c|c|c|c|} \hline
& \multicolumn{3}{r}{Reed Solomon $\times$ Goppa} & &
\multicolumn{3}{r}{Lomont Code \#1} & \\ \hline
rate & $k_1$ & $k_2$ & rate & $\delta$ & a & b & rate & $\delta$ \\
\hline \hline
0.1&81&  80&0.100046&0.470125&81 & 80&0.100099&0.469978\\
0.2&115&113&0.200633&0.306948&115&113&0.200015&0.307683\\
0.3&140&139&0.300448&0.205960&141&139&0.301183&0.205325\\
0.4&162&160&0.400185&0.136421&162&161&0.400443&0.136263\\
0.5&181&179&0.500216&0.086846&182&180&0.502632&0.085832\\
0.6&198&197&0.602223&0.051042&199&197&0.601205&0.051331\\
0.7&214&212&0.700448&0.027235&215&213&0.702037&0.026917\\
0.8&229&227&0.802578&0.011255&229&228&0.800183&0.011536\\
0.9&242&241&0.900448&0.002810&243&242&0.901015&0.002777\\
\hline
\end{tabular}
\end{table}

%\subsubsection{Optimal Rate Comparison} ??
\subsubsection{Optimal Rate Comparison}\label{s:CompareElliptic}
To show $err(q,\delta)<0$ for
all $q>1$ and $\delta\in[0,1]$, notice that
$$
\begin{array}{rcl}
q^2+q & > & q^2+q-2\\
q(q+1) & > & (q+2)(q-1)\\
\frac{q}{q-1} & > & \frac{q+2}{q+1}\\
\sqrt{\frac{q}{q-1}} & > & \sqrt{\frac{q+2}{q+1}}\\
\end{array}
$$
Thus the coefficient of $\delta$ in $err(q,\delta)$ \ref{f:err} is
always positive, so to find the maximum of $err(q,\delta)$, we can
assume $\delta=1$. Assuming $q>1$, if $err(q,1)<0$ then
$$
\begin{array}{rcl}
-\frac{2}{q^2-1}+2\left(\sqrt{\frac{q}{q-1}}-\sqrt{\frac{q+2}{q+1}}
\; \right) & < & 0
\end{array}
$$
Multiply each side by $\sqrt{(q-1)/q}(q+1)/2>0$, obtaining
$$
\begin{array}{rcl}
-\frac{1}{\sqrt{(q-1)q}} +1 + q -\sqrt{\frac{(q-1)(q+1)(q+2)}{q}} & < & 0\\
\left(1+q-\frac{1}{\sqrt{(q-1)q}}\right)^2 & < & \left(\sqrt{\frac{(q-1)(q+1)(q+2)}{q}}\right)^2\\
1-q-q^2+q^3+q^4-2\sqrt{(q-1)q}(q+1) & < & 2-q-3q^2+q^3+q^4\\
(2q^2-1)^2 & < & \left(2\sqrt{(q-1)q}(q+1)\right)^2\\
1+4q & < & 4q^3
\end{array}
$$
which always holds for $q>1$. The steps are reversible, proving
the optimal rate of the product code is slightly larger than the
optimal rate of the Lomont code \#2.

%\todo{can vary product code slightly to increase length?!
%nonrational point for divisor}

\subsubsection{The Degree 1 Case} I have been unable to complete
the analysis in this case. The trouble is computing
$k=h^0(C,S^a(\sheaf{E})\otimes\sheaf{O}(bP_0))$, for $\sheaf{E}$ a
rank 2, degree 1, indecomposable vector bundle on the elliptic
curve $C$. The problem is decomposing $S^a(\sheaf{E})$ in a manner
similar to theorem \ref{decomptheorem}. From the methods in the
classification this should be possible, but I have been unable to
solve it.

% end of file

\section{Conclusion} The codes in this chapter are comparable in performance
to product codes. Over the elliptic curves, since the degree 0
case was comparable to the product codes, perhaps the degree 1
case will be as good or better than the product codes. It would be
an interesting and worthwhile problem to finish this
classification, and see if there is any improvement.

\chapter{NEW CURVES OVER $P^2$}
%%%%%%%%%%%%%%%%%%%%%%%%%%%%%%%%%%%%%%%%%%%%%%%%%%%%%%%%%%%%%%%%%%%%%
% a paper on curves, from a computer search
% Chris Lomont Dec 2000
% accepted for publication Apr 17, 2002 - revisions needed - todo
%%%%%%%%%%%%%%%%%%%%%%%%%%%%%%%%%%%%%%%%%%%%%%%%%%%%%%%%%%%%%%%%%%%%%

%\NeedsTeXFormat{LaTeX2e}% LaTeX 2.09 can't be used (nor non-LaTeX)
%[1994/12/01]% LaTeX date must December 1994 or later

%\documentclass[12pt]{amsart}
%\usepackage{amscd}
%\usepackage{amsthm} % for proof environment
%\usepackage{amssymb}
%\usepackage{amsmath}
%\usepackage{amsfonts}

%%%%%%%%%%%%%%%%%%%%%%%%%%%%%%%%%%%%%%%%%%%%%%%%%%%%%%%%%%%%%%%%%%%%%
% place new commands here
%%%%%%%%%%%%%%%%%%%%%%%%%%%%%%%%%%%%%%%%%%%%%%%%%%%%%%%%%%%%%%%%%%%%%
%\newcommand{\field}[1]{\ensuremath{\mathbb{F}_{#1}}} %finite field

%%%%%%%%%%%%%%%%%%%%%%%%%%%%%%%%%%%%%%%%%%%%%%%%%%%%%%%%%%%%%%%%%%%%%
% start of document
%%%%%%%%%%%%%%%%%%%%%%%%%%%%%%%%%%%%%%%%%%%%%%%%%%%%%%%%%%%%%%%%%%%%%
%\begin{document}
\label{c:newcurves} This chapter is a paper that will be published
in \emph{Experimental Mathematics}.

% todo - table of contents pukes over this, so replace...
%\newcommand{\field}[1]{\ensuremath{\mathbb{F}_{#1}}} % finite field
\newcommand{\chfield}[1]{#1 Element Field} % finite field

{\center{\emph{\textbf{Yet More Projective Curves Over
\field{2}}}}}

\emph{\textbf{Abstract}} All plane curves of degree less than 7
with coefficients in \field{2} are examined for curves with a
large number of \field{q} rational points on their smooth model,
for $q=2^m, m = 3,4,\dots,11$. Known lower bounds are improved,
and new curves are found meeting or close to Serre's, Lauter's,
and Ihara's upper bounds for the maximal number of \field{q}
rational points on a curve of genus g.

\section{Introduction}

Let \field{q} denote the finite field with q elements. All
absolutely irreducible homogeneous polynomials $f \in
\field{2}[x,y,z]$ of degree less than 7 are examined for those
with a large number of \field{q} rational points, $q=2^m,
m=3,4,5,\dots,11$, extending the results in \cite{Moreno95}. A
brute force search obtained all rational points for each
polynomial of a given degree. The resulting list of polynomials
with many rational points, perhaps with singularities, were then
studied to determine if resolving singularities would add more
rational points on the smooth model. The result is an exhaustive
search of all curves resulting from desingularizing a homogeneous
polynomial of degree less than 7 in $\field{2}[x,y,z]$.

The rest of the paper is laid out as follows: first known bounds
on the maximal number of \field{q} rational points of a genus $g$
curve are recalled, along with some theorems that speed up the
computations. Then the computation is described in some detail. A
listing of the best found polynomials is given for each genus,
allowing checking (by computer unless one has a lot of time!) the
claimed number of rational points on each curve. Finally the new
lower bounds are listed in a table for each \field{q} and genus
combination.

\section{Genus Bounds and Irreducibility Tests}

Let $f \in \field{2}[x,y,z]$ be an absolutely irreducible
homogeneous polynomial; $f$ defines a projective plane curve $C$.
Let $\tilde{C}$ be the smooth model, and $g$ its genus. Some
bounds on the genus can be deduced from knowing the number of
\field{q} rational points in the plane and the number of
singularities in the plane. $N_q(g)$ is the maximum number of
\field{q}-rational points on a smooth curve of genus $g$ over
\field{q}. Serre's bound \cite{Serre} on $N_q(g)$ is
$$|N_q(g) - (q+1)|\leq g \lfloor 2\sqrt{q} \rfloor$$ where $\lfloor
\alpha \rfloor$ is the integral part of $\alpha$. This gives
$$N_q(g) \leq (q+1) + g\lfloor 2\sqrt{q}\rfloor$$ so if
there is an integer $g_0$ such that $$q+1+g_0\lfloor 2 \sqrt{q}
\rfloor < p$$ where $p$ is the point count on the particular curve
$C$ in question, then $g_0<g$. If the number of singularities is
$r$, and the degree of the polynomial $f$ is $d$, then
$$g\leq \frac{(d-1)(d-2)}{2}-r$$

To get an estimate of the total number of points possible on the
smooth model $\tilde{C}$ resulting from blowing up singularities
the following estimate was used.

\begin{thrm} Let $C\subseteq \mathbb{P}^2$
be a plane curve of degree d with singularities $P_1, P_2,
\dots,P_r$, with multiplicities $m_1,m_2,\dots,m_r$, for $r\geq
2$. \textbf{Then} $\sum^r_{i=1}m_i \leq \lfloor \frac{d}{2}
\rfloor r + 1$ if d is odd, and $\sum^r_{i=1}m_i \leq \lfloor
\frac{d}{2} \rfloor r$ if d is even.
\end{thrm}

So the number of points obtained from blowing up singularities is bounded above by $\lfloor \frac{d}{2} \rfloor r + 1$.

\begin{proof}
By Bezout's theorem, a line through any 2 singularities $P_i$ and
$P_j$ implies $m_i + m_j \leq d$. Thus at most one singularity can
have multiplicity $> d/2$, and the result follows.
\end{proof}

More details on resolution of curve singularities can be found in
\cite{Hartshorne} and \cite{Walker}.

To test for absolute
irreducibility the following was used \cite{Ragot}:

\begin{defn}
Let $k$ be a field. The polynomial $f\in k[x_1,x_2,\dots,x_n]$ has
a {\bf simple solution} at a point $P\in k^n$ if $f\in
I(P)\diagdown I(P)^2$, with $I(P)$ being the ideal of polynomials
vanishing at $P$.\end{defn}

\begin{thrm} \textbf{If} $f\in k[x_1,x_2,\dots,x_n]$ is irreducible over the
perfect field $k$ and has a simple solution in $k^n$,
\textbf{then} $f$ is absolutely irreducible.\end{thrm}
\begin{proof} Since $f$ is irreducible over $k$, its absolutely
irreducible factors are conjugate over $k$. If $P\in k^n$ is a
root of one of the factors, it must be a root of the others. But
$P$ only vanishes to order 1, thus $f$ has one factor, and is
absolutely irreducible. \end{proof}

The number of \field{q} points on a plane curve can be computed
directly by brute force, as can lower bounds on the number of
singularities, and then the above inequalities can be used to
obtain bounds on the genus. So given a polynomial $f$, the number
of rational points computed on it, and the number of singularities
found, upper and lower bounds on the possible genus are obtained,
which speeds up the search by removing curves early in the
computation that have uninteresting combinations of genus and
rational point count.

\section{Computation}
\subsection{Storing Polynomials Compactly}
All \field{q} rational points were found on each homogeneous
polynomial of degree $\leq 5$ in $\field{2}[x,y,z]$, for $q=2^m,
m=3,4,\dots,11$. Due to the the time required, degree 6
homogeneous polynomials were examined only for $m=3,4,\dots,9$.

The most time consuming part was counting the number of
\field{q}-rational points on each plane curve. This was done with
a C program using exhaustive search. Several ideas were used to
reduce the complexity at each stage. Degree 6 computations will be
described; the other degrees are similar. The code was checked for
correctness by comparing the degree $\leq$ 5 results with
\cite{Moreno95}, and in the process a few curves were found that
were previously overlooked.

First, each homogeneous polynomial can be represented uniquely by
a 32 bit integer, using each bit to signify the presence of a
certain monomial in the polynomial. In degree 6, there are
$\binom{6+2}{2}=28$ different monomials of the form $x^i y^j z^k$
with $i+j+k=6$ and $0\leq i,j,k$. Each bit from 0 to 27 denotes
the presence of a monomial, and the mapping
\makebox{$\alpha:\{homogeneous \,\,degree \,\, 6 \,\,f\in
\field{2}[x,y,z]\} \rightarrow \{1,2,\dots,2^{28}-1 \}$} thus
defined is a bijection. Thus each homogeneous polynomial
\makebox{$f\in \field{2}[x,y,z]$} of degree 6 corresponds to a
unique integer between $1$ and $2^{28}-1\approx268$ million.

\subsection{Reducing Computation Time}
To reduce the number of polynomials searched, equivalent ones
under the action of $GL_3(\field{2})$ on the variables $x,y,z$
were removed. To fit the entire degree 6 computation in memory, a
bit table of 32 megabytes of RAM was used, with the position of
each bit representing the number of a polynomial using the above
bijection. All bits were set to 1, denoting all polynomials are
still in the search space. Then the orbit of each polynomial under
$GL_3(\field{2})$ was removed from the bit table, and the
polynomial in each orbit requiring the least computation to
evaluate was written to a data file. Since the size of
$GL_3(\field{2})$ is 168, this was expected to give approximately
a 168 fold decrease in the number of curves needing to be searched
(not exactly 168 since some polynomials are invariant under some
automorphisms). By using the representative of each orbit
requiring the least work to evaluate, the search time was reduced
significantly (see below). This trimmed the 268 million degree 6
polynomials down to 1.6 million. Also, clearly reducible
polynomials, such as those with all even exponents or divisible by
a variable, were removed at this point. At each stage data was
saved to prevent having to rerun any step.

For speed reasons finding solutions was done by table lookup, so
in each orbit the polynomial needing the fewest number of lookups
was selected. By choosing the representative with the fewest
number of lookups as opposed, for example, to the polynomial with
the lowest value of $\alpha(f)$ defined above, 12 million lookups
were removed from polynomials of degree 6, resulting in over 3
trillion operations removed during rational point counting.

\subsection{Timing}
After the C program computes all the \field{q}-rational points,
the points are tested for singularities (a singularity will add an
additional \field{q} rational point only if it comes from
resolving a \field{q} rational singularity). The computation up to
this point took about 80 hours of computer time on a Pentium III
800 MHZ. Using the bounds above on the genus and possible ranges
for number of \field{q} rational points on the smooth model, the
program searched all curves for those with a large number of
possible \field{q} rational points for each genus and field
combination, and all such curves were written out to be examined.
If the genus of one of these curves was not forced to be unique
using the bounds, the program KANT \cite{Kant} was used to compute
the genus, and this data was incorporated into the C program, and
another pass was run. Due to the large number of degree 6 curves,
and the length of time to compute the genus of all of them, not
all degree 6 curves of genus $\leq 5$ were identified. All curves
of degree 6, genus $\geq 6$ were identified. The C program also
found simple points over \field{2} to apply the irreducibility
theorem above, and then Maple V \cite{Maple} was used to test for
irreducibility since it has multivariable factoring algorithms
over finite fields. For 12 curves of degree 6, there were no
simple \field{2} points, so \field{4} simple points were used. For
2 of these curves there were no such simple \field{4} points, so
\field{8} simple points were used. This turned out to suffice to
check absolute irreducibility of all polynomials in this paper.
The C program also found the singularity types of the \field{2}
singularities for visual inspection to see if there were clearly
more rational points on the smooth model. The package of
\cite{Hache} was not available to do a more detailed singularity
analysis, thus some of the bounds below may be improved by looking
for rational points over a wider class of singularities than the
\field{2} singularities considered here.

The final C code can be found at \cite{LomontSource}.

\section{Computational Results}

For each field and genus combination polynomials are listed that
result in the largest found number of rational points on the
smooth model of the curve. For fields \field{q}, $q=2^m,
m=3,4,\dots,11$, all homogeneous polynomials in $\field{2}[x,y,z]$
of degree $\leq 5$ were searched. For $m=3,4,\dots,9$, the search
was extended to include all degree 6 homogeneous polynomials in
$\field{2}[x,y,z]$. For genus and field combinations not listed
here, see \cite{Moreno95}.

\emph{Remark:} the four polynomials found in \cite{Moreno95} of
degree 4, genus 3, with 113 rational points over \field{64}, are
only 2 distinct polynomials modulo the action of $GL_3(\field{2})$
on the variables $x,y,z$.

\subsection*{\chfield{8}}
A curve of genus 3 and the maximal number of smooth points, 24, is
the Klein Quartic $$x^3\,y + y^3\,z + x\,z^3$$
%%$$ x^6 + x^5\,y + x\,y^5 + y^6 +
%%  \left( x^5 + x^2\,y^3 + y^5 \right) \,z + x^3\,y\,z^2 +
%%  x\,y^2\,z^3 + \left( x + y \right) \,z^5 + z^6$$
A genus 5 curve with 28 planar smooth points is $$x^6 + x^5\,y + x^3\,y^3 + y^6 + y^5\,z +
y^4\,z^2 + \left( x^3 + x\,y^2 + y^3 \right) \,z^3 +
  \left( x^2 + x\,y \right) \,z^4 + x\,z^5$$
\emph{Note}: a reviewer remarked that a genus 5 curve is known
  with 29 points \cite{vdvg}.\\

\noindent A genus 6 curve, with 33 planar smooth points is
$$x^4\,y^2 + x^3\,y^3 + x\,y^5 + x\,y^4\,z + y^4\,z^2 +  \left(
x^2 + y^2 \right) \,z^4 + y\,z^5$$ A curve of genus 7 with 33
smooth planar points is
  $$x^6 + x^5\,y + x^4\,y^2 + x^3\,y^3 + y^6 + y^5\,z +
  \left( x^2\,y^2 + x\,y^3 + y^4 \right) \,z^2 +
  \left( x^2\,y + x\,y^2 \right) \,z^3 +
  \left( x^2 + x\,y + y^2 \right) \,z^4$$
A genus 8 curve with 33 smooth planar points is $$x^5\,y
+ x^2\,y^4 + x\,y^5 + y^6 +
  \left( x^4\,y + x^2\,y^3 \right) \,z + x^3\,y\,z^2 +
  \left( x^3 + x\,y^2 \right) \,z^3 + x\,y\,z^4 + y\,z^5$$
Two curves of genus 9, each with 33 smooth planar points:
$$f_1=x^5\,y + x^4\,y^2 + x^2\,y^4 +
  \left( x^3\,y^2 + x^2\,y^3 \right) \,z + x^4\,z^2 +
  \left( x\,y^2 + y^3 \right) \,z^3 + x^2\,z^4 + x\,z^5 +
  z^6$$ $$f_2=x^6 + x^4\,y^2 + x^3\,y^3 + x^2\,y^4 + x^3\,y^2\,z +
  \left( x^4 + x^3\,y + x\,y^3 \right) \,z^2 + y^3\,z^3 +
  \left( x + y \right) \,z^5 + z^6$$
Five curves of genus 10 with 35 smooth points in the plane:
  $$f_1 = x^5\,y + x^2\,y^4 + y^6 +
  \left( x^3\,y^2 + x\,y^4 + y^5 \right) \,z + y^4\,z^2 +
  \left( x^2\,y + x\,y^2 \right) \,z^3 +
  \left( x^2 + y^2 \right) \,z^4 + x\,z^5$$ $$f_2 = x^5\,y + x^4\,y^2 +
  x^3\,y^3 + x\,y^5 + y^6 + x^2\,y^3\,z + \left( x^4 + x^2\,y^2 + x\,y^3 \right) \,
   z^2 + x^3\,z^3 + y^2\,z^4 + x\,z^5 + z^6$$ $$f_3 = x^5\,y + x^4\,y^2 + x^2\,y^4 +
  \left( x^4\,y + y^5 \right) \,z +
  \left( x^4 + x\,y^3 \right) \,z^2 + x^3\,z^3 +
  \left( x^2 + y^2 \right) \,z^4 + y\,z^5 + z^6 $$ $$f_4 = x^5\,y + x^3\,y^3 + x^2\,y^4 +
  \left( x^5 + x^2\,y^3 + y^5 \right) \,z +
  \left( x^2\,y + y^3 \right) \,z^3 + x^2\,z^4 +
  \left( x + y \right) \,z^5 + z^6 $$ $$f_5 = x^4\,y^2 + x^2\,y^4 + x\,y^5 +
  x^5\,z + x\,y^3\,z^2 + \left( x^3 + x^2\,y \right) \,z^3 +
  \left( x^2 + y^2 \right) \,z^4 + x\,z^5 $$

\subsection*{\chfield{16}}
One genus 6 curve, a Hermitian curve, with the maximal number of
smooth points, 65, was found (it is known to be the unique such
curve up to isomorphism):
$$x^5 + y^5 + z^5$$ Two genus 7 curves each with 57 smooth points
in the plane:
\begin{eqnarray*} f_1&=&x^4\,y^2 + x\,y^5 + y^6 +
  \left( x^2\,y^3 + x\,y^4 + y^5 \right) \,z +
  \left( x^2 + x\,y \right) \,z^4 + y\,z^5 \\  f_2&=&x^4\,y^2 + x^2\,y^4 + y^6 +
  \left( x^2\,y^3 + x\,y^4 \right) \,z +
  \left( x^2 + x\,y \right) \,z^4 + y\,z^5\end{eqnarray*}
A curve with genus 8 with 57 smooth plane points is $$x^6 +
x^3\,y^3 + x^2\,y^4 + x^4\,y\,z + x^2\,y^2\,z^2 + \left( x^3 +
x^2\,y \right) \,z^3 +
  \left( x + y \right) \,z^5$$
There are two curves of genus 9 with 57 smooth plane points, each
receiving two points from blowups: $f_1$ from the singularity
$(1:1:1)$ of type $u^2+u v + v^2$ which splits over \field{16},
and $f_2$ from the singularity $(0:1:1)$ of type $u v$. Thus
$N_{16}(9)\geq 59$.
\begin{eqnarray*} f_1&=&x^5\,y + x^3\,y^3 + x\,y^5 +
  \left( x^5 + y^5 \right) \,z + x^2\,y^2\,z^2 +
  \left( x^3 + y^3 \right) \,z^3 +
  \left( x + y \right) \,z^5 \\ f_2 &=& x^6 + x^5\,y + x^2\,y^4 + y^5\,z + x^2\,y^2\,z^2 +
  x\,y^2\,z^3 + \left( x^2 + x\,y \right) \,z^4 + y\,z^5\end{eqnarray*}
The two curves of genus 10 each with 59 plane smooth points are:
  $$f_1 = x^5\,y + y^6 + \left( x^2\,y^3 + y^5 \right) \,z +
  \left( x^4 + x^3\,y + x\,y^3 \right) \,z^2 +
  x\,y^2\,z^3 + \left( x + y \right) \,z^5 + z^6$$
  $$f_2 = x^5\,y + y^6 + \left( x^4\,y + x\,y^4 + y^5 \right) \,z +
  \left( x^4 + x\,y^3 \right) \,z^2 +
  \left( x^3 + x\,y^2 \right) \,z^3 + y^2\,z^4 + z^6$$

\subsection*{\chfield{32}}
Three curves with genus 4 and 71 smooth points on the plane curve
are: \begin{eqnarray*} f_1 &=& x^4\,y + x\,y^4 + y^5 + x\,y^3\,z +
  \left( x\,y^2 + y^3 \right) \,z^2 + x^2\,z^3 + x\,z^4 +
  z^5 \\ f_2&=& x^6 + x^3\,y^3 + y^6 + \left( x^4\,y + y^5 \right) \,z +
  \left( x^3\,y + x^2\,y^2 \right) \,z^2 +
  \left( x^3 + x^2\,y + y^3 \right) \,z^3 + \\&& x^2\,z^4 +
  y\,z^5 + z^6 \\ f_3 &=& x^6 + x^3\,y^3 + y^6 + \left( x^5 + x^3\,y^2 +
     x^2\,y^3 \right) \,z + y^4\,z^2 +
  \left( x^3 + y^3 \right) \,z^3 + y^2\,z^4 + x\,z^5 + z^6
  \end{eqnarray*}
A curve with 82 smooth points in the plane and genus 5 is
  $$x^6 + x^3\,y^3 + x^2\,y^4 + y^6 + x^5\,z + x^3\,y\,z^2 +
  \left( x^3 + x\,y^2 + y^3 \right) \,z^3 + x^2\,z^4 +
  y\,z^5$$
A genus 6 curve with 82 planar smooth points and 2 points above
the singularity $(1:0:1)$ of type $u v$ (thus $N_{32}(6)\geq 84$)
is $$x^6 + y^6 + \left( x^4\,y + x^3\,y^2 + x\,y^4 \right) \,
   z + x\,y^2\,z^3 + \left( x^2 + x\,y + y^2 \right) \,z^4$$
Two genus 7 curves each with 92 planar smooth points are
\begin{eqnarray*}f_1&=&x^3\,y^3 + y^6 + \left( x^5 + x^3\,y^2 \right) \,z +
  \left( x^4 + y^4 \right) \,z^2 +
  \left( x^3 + y^3 \right) \,z^3 + y^2\,z^4 + x\,z^5 + z^6\\
   f_2&=&x^6 + y^6 + \left( x^5 + y^5 \right) \,z + y^4\,z^2 +
  \left( x^3 + y^3 \right) \,z^3 + x^2\,z^4 +
  \left( x + y \right) \,z^5\end{eqnarray*}
A curve with 93 planar smooth points, genus 8, is $$x\,y^5 + y^6 + \left( x^5 + x^4\,y
\right) \,z + y^4\,z^2 + \left( x^3 + y^3 \right) \,z^3 + y^2\,z^4
+ y\,z^5$$
A genus 9 curve with 93 smooth planar points:  $$x^4\,y^2 +
x^3\,y^3 + \left( x^5 + x^3\,y^2 + x\,y^4 +
     y^5 \right) \,z + x^2\,y^2\,z^2 +
  \left( x^3 + y^3 \right) \,z^3 + x^2\,z^4 + z^6$$
Genus 10 with 103 smooth planar points:
  $$x^6 + x^3\,y^3 + x\,y^5 +
  \left( x^2\,y^2 + x\,y^3 \right) \,z^2 +
  \left( x^3 + x\,y^2 + y^3 \right) \,z^3 + x\,y\,z^4 +
  \left( x + y \right) \,z^5$$

\subsection*{\chfield{64}}
One curve had genus 4 and 118 smooth planar points: $$x^3\,y^2 +
y^5 + y^4\,z + y^2\,z^3 + z^5$$ Two curves of genus 6 had 160
smooth planar points (which is one less than the bound of 161):
\begin{eqnarray*} f_1 &=& x^4\,y^2 + x^2\,y^4 + x\,y^5 + y^5\,z +
y^3\,z^3 +
  y\,z^5 + z^6 \\ f_2&=&x^6 + x^5\,z + \left( x^4 + y^4 \right) \,z^2 + x^3\,z^3 +
  y^2\,z^4 + y\,z^5\end{eqnarray*}
Genus 7, 153 planar smooth points: $$x^2\,y^4 + x\,y^5 + y^6 +
  \left( x^3\,y^2 + x\,y^4 + y^5 \right) \,z +
  x\,y^3\,z^2 + x^2\,z^4 + x\,z^5 + z^6$$
Three curves had genus 8 and 159 plane smooth points, the last two
of which have no rational points over \field{2}:
\begin{eqnarray*} f_1 &=& x^3\,y^3 + y^6 + \left( x^4\,y + x\,y^4
\right) \,z +
  \left( x^3 + y^3 \right) \,z^3 + x\,y\,z^4 \\ f_2 &=& x^6 + x^5\,y + x^3\,y^3 + x\,y^5 + y^6 +
  \left( x^3\,y^2 + y^5 \right) \,z + x^3\,y\,z^2 +
  \left( x^2\,y + y^3 \right) \,z^3 + \\&& y^2\,z^4 + x\,z^5 +
  z^6 \\ f_3 &=& x^6 + x^4\,y^2 + x^3\,y^3 + x^2\,y^4 + y^6 +
  \left( x^4 + x^2\,y^2 + y^4 \right) \,z^2 +
  \left( x^3 + y^3 \right) \,z^3 + \\&&
  \left( x^2 + y^2 \right) \,z^4 + z^6\end{eqnarray*}
There are 166 plane smooth points on this curve of genus 9: $$x^6
+ x^3\,y^3 + \left( x^4\,y + x^2\,y^3 \right) \,z +
  \left( x^3\,y + x\,y^3 + y^4 \right) \,z^2 + x^2\,z^4 +
  y\,z^5$$
Four genus 10 curves each had 171 points on their smooth model:
\begin{eqnarray*} f_1 &=& x^6 + y^6 + \left( x^4\,y + x^2\,y^3 + x\,y^4 \right)
\,
   z + x^3\,y\,z^2 + x\,y^2\,z^3 + x\,y\,z^4 + z^6 \\ f_2 &=& x^6 + x^5\,y + x^4\,y^2 + x^3\,y^3 + x^2\,y^4 + x\,y^5 +
  y^6 + \left( x^4\,y + x\,y^4 \right) \,z +
  \left( x^2\,y + x\,y^2 \right) \,z^3 + z^6\\ f_3 &=& x^6 + x^3\,y^3 + y^6 + \left( x^4\,y + x\,y^4 \right) \,
   z + \left( x^3 + y^3 \right) \,z^3 + z^6 \\f_4 &=& x^6 + x^3\,y^3 + x\,y^5 + x^3\,y^2\,z +
  \left( x^4 + x^3\,y + y^4 \right) \,z^2 + y^2\,z^4 +
  x\,z^5 + z^6 \end{eqnarray*}

\subsection*{\chfield{128}}
There is one degree 6 plane curve with a genus 3 smooth model,
with 183 smooth plane points, and another point coming from the
singularity $(0:0:1)$ of type $(u+v)(u^2+u v + v^2)$, which
matches \cite{Moreno95}. The curve is $$x^6 + x^5\,y + x^4\,y^2 +
x^3\,y^3 + x^2\,y^4 +
  \left( x^5 + x^4\,y \right) \,z + y^4\,z^2 +
  \left( x^3 + y^3 \right) \,z^3$$
A curve of genus 4 with 215 planar smooth points (2 less than the maximum possible) is $$x^2\,y^3 +
x\,y^4 + x^4\,z + x\,y^2\,z^2 + x\,y\,z^3 +
  \left( x + y \right) \,z^4$$
There are two curves of genus 6 with 240 planar smooth points,
receiving 3 points each from singularities. $f_1$ has type $u
v(u+v)$ at $(0:1:0)$ and $f_2$ has type $(u+v)(u^2+u v+v^2)$ at
$(0:1:0)$ and type $u v$ at $(1:0:0)$. Thus $N_{128}(6)\geq 243$.
\begin{eqnarray*} f_1 &=& x^4\,y^2 + \left( x^5 + x^4\,y + x^2\,y^3 \right) \,z +
  \left( x^2\,y^2 + x\,y^3 \right) \,z^2 +
  \left( x^2 + x\,y + y^2 \right) \,z^4 \\ f_2 &=& x^3\,y^3 + x^4\,y\,z + \left( x^4 + x^3\,y \right) \,z^2 +
  \left( x^3 + x^2\,y + y^3 \right) \,z^3 + z^6\end{eqnarray*}
Two genus 7 curves with 248 smooth planar points:
\begin{eqnarray*} f_1 &=& x^3\,y^3 + x\,y^5 + y^6 + x^3\,y^2\,z + y^4\,z^2 +
  x^3\,z^3 + \left( x^2 + y^2 \right) \,z^4 +
  \left( x + y \right) \,z^5 \\ f_2 &=& x^5\,y + x^4\,y^2 + x^2\,y^4 +
  \left( x^3\,y^2 + x^2\,y^3 \right) \,z + x^4\,z^2 +
  x\,y^2\,z^3 + z^6 \end{eqnarray*}
A curve with 266 planar smooth points, genus 8, and no \field{2}
rational points is $$x^6 + x^3\,y^3 + x^2\,y^4 + x\,y^5 + y^6 +
  \left( x^5 + y^5 \right) \,z +
  \left( x^2\,y^2 + y^4 \right) \,z^2 +
  \left( x^3 + y^3 \right) \,z^3 + x\,z^5 + z^6$$
There are 269 smooth plane points on the curves of genus 9 given
by \begin{eqnarray*} f_1 &=& x^4\,y^2 + x\,y^5 + \left( x^4 + y^4
\right) \,z^2 +
  \left( x^3 + y^3 \right) \,z^3 + x\,y\,z^4 + x\,z^5 + z^6 \\ f_2
  &=& x^6 + x^3\,y^3 + x^2\,y^4 + y^6 +
  \left( x\,y^4 + y^5 \right) \,z + x^2\,y^2\,z^2 +
  \left( x^2 + x\,y + y^2 \right) \,z^4 \end{eqnarray*}
The smooth curve of genus 10 with 276 \field{128} rational points
is $$x^6 + y^6 + x^2\,y^3\,z +
  \left( x^4 + x^3\,y + y^4 \right) \,z^2 + x^3\,z^3 +
  x^2\,z^4 + x\,z^5$$

\subsection*{\chfield{256}}
A genus 3 curve not listed in \cite{Moreno95} with 350 smooth
planar points is given by $$x^5 + x\,y^4 + y^5 + \left( x^2\,y^2 +
y^4 \right) \,z +
  \left( x^2\,y + x\,y^2 \right) \,z^2 + x\,z^4 + z^5$$
A curve with 399
smooth plane points and genus 5 is $$x^6 + x^4\,y^2 + x^5\,z +
  \left( x^2\,y^2 + y^4 \right) \,z^2 +
  \left( x^2\,y + x\,y^2 \right) \,z^3 + x^2\,z^4 + y\,z^5$$
A genus 6 curve with 416 smooth plane points is $$x^4\,y +
x^3\,y^2 + y^4\,z +
  \left( x^3 + y^3 \right) \,z^2 +
  \left( x^2 + x\,y \right) \,z^3 + z^5$$
One point from the singularity $(1:0:0)$ of type $u v^2$ is
added to the 442 smooth plane points on a curve of genus 7 given
by $$x^3\,y^3 + x^2\,y^4 + y^5\,z + x^3\,y\,z^2 +
  \left( x\,y^2 + y^3 \right) \,z^3 + y^2\,z^4 + z^6$$
A curve of genus 8 with one point less than the Serre bound has
512 smooth plane points and is given by $$x^4\,y^2 + y^5\,z +
x\,z^5$$ Two curves of genus 9, each with 474 smooth points and 2
points from singularities of type $u^2+u v + v^2$, which factor
over \field{256}, at points $(0:1:1)$ and $(1:1:0)$ respectively
(so $N_{256}(9) \geq 476$) are \begin{eqnarray*}f_1&=&x^5\,y +
x^3\,y^3 + x^2\,y^4 + x\,y^5 + y^4\,z^2 + x^3\,z^3 + y^2\,z^4 +
x\,z^5 \\f_2&=&x^6 + y^6 + \left( x^5 + y^5 \right) \,z + x^4\,z^2
+ \left( x^2\,y + x\,y^2 \right) \,z^3 + x\,z^5 +
z^6\end{eqnarray*} Two smooth curves of genus 10 have 537 smooth
plane points: \begin{eqnarray*} f_1&=&x^6 + x\,y^5 + x^4\,y\,z +
x^2\,y^2\,z^2 + y^3\,z^3 +
  x\,z^5 \\f_2&=&x^6 + x^5\,y + x^3\,y^3 + x\,y^5 + y^6 +
  \left( x^5 + y^5 \right) \,z + x^2\,y^2\,z^2 +
  x\,y\,z^4 + z^6\end{eqnarray*}
\subsection*{\chfield{512}}
Four curves overlooked in \cite{Moreno95} of genus 4 have 663
plane smooth points. They are
\begin{eqnarray*}f_1 &=& x^4\,y + x\,y^4 + \left( x^3\,y + y^4 \right) \,z +
  \left( x\,y + y^2 \right) \,z^3 + z^5 \\ f_2&=&x^4\,y + x\,y^4 + y^5 +
  \left( x\,y^3 + y^4 \right) \,z + \left( x\,y^2 + y^3 \right) \,z^2 +
  x^2\,z^3 + x\,z^4 \\ f_3&=&x^4\,y + x^2\,y^3 + y^5 +
  \left( x^2\,y^2 + x\,y^3 + y^4 \right) \,z +  \left( x^3 + y^3 \right)
  \,z^2 + z^5 \\ f_4&=&x^5 + y^5 + \left( x^4 + x^3\,y + y^4 \right) \,z +
  \left( x\,y^2 + y^3 \right) \,z^2 + z^5\end{eqnarray*}
A genus 6 curve with 766 smooth plane points and one more point
from the singularity $(1:0:0)$ of type $(u+v)(u^2+u v+v^2)$ (so
$N_{512}(6)\geq 767$) is $$x^3\,y^3 + y^6 + \left( x\,y^4 + y^5
\right) \,z +
  \left( x^2\,y^2 + y^4 \right) \,z^2 + x^3\,z^3 +
  x\,y\,z^4 + \left( x + y \right) \,z^5$$
There are 786 smooth plane points and 1 point from the
singularity $(1:0:0)$ of type $u v^2$ on the genus 7 curve
$$x^2\,y^4 + y^6 + x^3\,y^2\,z +
  \left( x^3 + x\,y^2 \right) \,z^3 + x\,y\,z^4 + y\,z^5$$
A curve of genus 8 with 813 plane smooth points is $$x^2\,y^4 +
y^6 + \left( x^5 + x^2\,y^3 \right) \,z +
  \left( x^3\,y + x\,y^3 + y^4 \right) \,z^2 +
  \left( x^3 + x^2\,y \right) \,z^3 + x\,z^5$$
A genus 9 curve with 837 smooth plane points is $$x^6 + x^4\,y^2 +
\left( x\,y^4 + y^5 \right) \,z +
  x^2\,y^2\,z^2 + \left( x^2\,y + x\,y^2 \right) \,z^3 +
  x\,y\,z^4 + x\,z^5 + z^6$$
A smooth genus 10 plane curve with 845 plane points is $$x^5\,y +
x^4\,y^2 + x^2\,y^4 + y^6 +
  \left( x^2\,y^3 + y^5 \right) \,z +
  \left( x^3\,y + y^4 \right) \,z^2 + x^2\,z^4 +
  \left( x + y \right) \,z^5$$

\subsection*{\chfield{1024}}
A genus 3 curve with 1211 smooth plane points is $$x^3\,y + y^3\,z
+ y\,z^3 + z^4$$ Three genus 4 curves have 1273 smooth plane
points:
\begin{eqnarray*} f_1 &=& x^3\,y^2 + y^5 + x^2\,y\,z^2 + y^2\,z^3 + x\,z^4\\
f_2 &=& x^4\,y + x^2\,y^3 + y^5 +
  \left( x^2\,y^2 + y^4 \right) \,z +
  \left( x^3 + y^3 \right) \,z^2 + x\,z^4 \\
  f_3 &=& x^5 + x^3\,y^2 + x^2\,y^3 + y^5 + y^4\,z + x\,y^2\,z^2 +
  x^2\,z^3\end{eqnarray*}
A curve with 1343 smooth plane points, genus 5, and 2 points
coming from the singularity $(1:0:0)$ of type $u v$ (and thus attaining the maximum possible 1345) is
$$x^3\,y^2 + y^5 + x^3\,y\,z + y^3\,z^2 + z^5$$
A genus 6 curve
with 1383 smooth plane points is $$x^4\,y + x\,y^4 + y^5 +
x^2\,y^2\,z + x\,y\,z^3 + z^5$$

\subsection*{\chfield{2048}}
Two genus 3 curves with 2293 smooth plane points, and one more
coming from the singularity $(0:1:0)$ of type $(u+v)(u^2+u v+v^2)$
on each curve are $$f_1=x^4\,y + x^3\,y^2 + x^3\,y\,z + y^2\,z^3 +
  \left( x + y \right) \,z^4$$ $$f_2=x^4\,y + x^3\,y^2 + x^3\,y\,z + x^3\,z^2 + y^2\,z^3 +
  y\,z^4$$
Three curves with 2380 smooth plane points and genus 4 are
\begin{eqnarray*}f_1&=&x^5 + y^5 + x^3\,y\,z + y^2\,z^3 + x\,z^4 \\ f_2&=&x^5 +
x^4\,y + \left( x\,y^3 + y^4 \right) \,z + x^3\,z^2 + \left( x^2 +
y^2 \right) \,z^3 + z^5 \\ f_3&=&x^5 + x^2\,y^3 + x\,y^4 + x^4\,z
+ x^3\,z^2 +  \left( x^2 + x\,y + y^2 \right) \,z^3 + x\,z^4 + z^5
  \end{eqnarray*}
A genus 5 curve with 2422 smooth plane points is $$x^4\,y +
x^3\,y^2 + y^5 + x\,y^2\,z^2 +
  \left( x^2 + x\,y + y^2 \right) \,z^3$$
Finally, a genus 6 curve with 2556 planar smooth points is
$$x^4\,y + x^2\,y^3 + x\,y^4 + y^5 +
  \left( x^3\,y + y^4 \right) \,z + x^2\,z^3 + y\,z^4$$

\section{Tallies}

The columns headed ``bound" give the Serre \cite{Serre} bound,
unless marked \cite{Ihara} as Ihara or \cite{Lauter} as Lauter.
The columns headed ``best" give the lower bounds for $N_q(g)$
found above; bounds marked \cite{Moreno95} and \cite{Serre} are
from those previous papers. Note in particular the reduction from
the Serre \cite{Serre} upper bounds using \cite{Ihara} and
\cite{Lauter} has made several known curves closer to or already
optimal. After this paper was initially written in 2000, improved
bounds were published in \cite{HL} and incorporated in this table.
These new bounds made the $q=128$, $g=4$ curve optimal.

% set table spacing
 \setlength{\tabcolsep}{0.6mm}

\begin{table}
\caption{Known Bounds on Number of Rational Points on
Curves}\label{T:tallies}
\begin{tabular}{|c|c|c|c|c|c|c|c|c|} \hline
$F_q$ & best 3 & bound  3 & best 4 & bound 4 & best 5 & bound 5 &
best 6 & bound 6
\\ \hline

8 & 24\cite{Serre} & 24 & 25\cite{Moreno95} & 25\cite{HL} & 28 &
30 \cite{HL} & 33 & 35 \cite{Lauter}
\\ \hline

16 & 38\cite{Serre} & 41 & 45\cite{Moreno95} & 45\cite{HL} &
45\cite{Moreno95} & 53\cite{HL} & 65 & 65 \\ \hline

32 & 63\cite{Moreno95} & 65\cite{Lauter} & 71 & 74\cite{HL} & 82 &
85\cite{HL} & 84 & 96\cite{HL} \\ \hline

64 & 113\cite{Moreno95} & 113 & 118 & 129 & 130\cite{Moreno95} &
145 & 160 & 161
\\ \hline

128 & 184\cite{Moreno95} & 195 & 215 & 215 \cite{HL} &
227\cite{Moreno95} & 239 & 243 & 258 \cite{HL}
\\ \hline

256 & 350\cite{Moreno95} & 353 & 381\cite{Moreno95} & 385 & 399 &
417 & 416 & 449
\\ \hline

512 & 640\cite{Moreno95} & 648 & 663 & 693 & 724\cite{Moreno95} &
738 & 767 & 783 \\ \hline

1024 & 1211 & 1217 & 1273 & 1281 & 1345 & 1345 & 1383 & 1409 \\
\hline

2048 & 2294 & 2319 & 2380 & 2409 & 2422 & 2499 & 2556 & 2589 \\
\hline

\hline

$F_q$ & best 7 & bound  7 & best 8 & bound 8 & best 9 & bound 9 &
best 10 & bound 10
\\ \hline

8 & 33 & 38 \cite{HL} & 33 & 42\cite{HL} &33 & 45\cite{HL} &35 &
49\cite{HL} \\ \hline

16 & 57 & 69\cite{HL} &57 &75\cite{HL} &59 &81\cite{Ihara} &59
&87\cite{Ihara}
\\ \hline

32 & 92 & 107 \cite{HL} & 93 & 118 \cite{HL} & 93 & 128 \cite{HL} & 103 & 139 \cite{HL}\\
\hline

64 & 153 &177 & 159 &193 & 166 &209  & 171  & 225\\ \hline

128 & 248 & 283 & 266 &302\cite{HL} & 269 & 322 \cite{HL} & 276 & 349 \\
\hline

256 & 443 & 481 & 512 & 513 & 476 & 545 & 537 & 577 \\ \hline

512 & 787 & 828 & 813 & 873  & 837 & 918  & 845 & 963 \\ \hline

\end{tabular}
\end{table}

\section{Comments} The techniques used here make a search over
degree 7 plane curves feasible on a supercomputer, and quite
possibly on a home PC. The desingularized curves can be used to
construct algebraic-geometric Goppa codes \cite{Pretzel},
\cite{TV}. For example, using the genus 5 curve over \field{1024}
with 1345 rational points, linear codes with parameters
$[n,k-4,n-k]$ can be constructed for $10\leq n \leq 1344$ and
$8<k<n$ over \field{1024}. Similarly, using the genus 6 \field{64}
curve with 160 points, $[n,k-5, n-k]$-linear codes can be
constructed for $12\leq n \leq 159$ and $10< k <n$, and the
\field{256} curve of genus 8 with 512 rational points gives
$[n,k-7,n-k]$-linear codes for $16\leq n\leq 511$ and $14<k<n$
(see, for example, \cite{Pretzel}).

\emph{Thanks} to the reviewer for numerous suggestions on layout
and a few corrections.

For conclusions and open problems see section
\ref{s:NewCurveConclusion}.

% end of file

\chapter{CONCLUSION AND OPEN PROBLEMS}
% conclusion of Chris Lomont's Thesis 2003
\todo{make headings here consistent}
\section{Codes Obtained}
Codes resulting from ruled surfaces over $\proj^1$ are completely
classified in this thesis. Although these Lomont codes \#1 still
suffer from having a fixed length like the Reed Solomon codes, the
Lomont codes \#1 are usable in some situations where the product
Reed Solomon fails. For example, the Lomont codes \#1 can correct
longer burst errors over the same fixed field, perhaps making
hardware implementations more cost efficient. Decoding them using
the methods for a (product) of Goppa codes is still not as
efficient as decoding Reed Solomon codes, so there is still work
to do on decoding.

For the codes resulting from ruled surfaces over elliptic curves,
the degree 1 bundle case is left open, and is tantalizingly close
to being solvable. To do this one needs some classification of
symmetric powers similar to theorem \ref{decomptheorem}, but for
the rank 2 indecomposable bundle of degree 1. The codes from
degree 0 bundles are comparable to the product code of a Reed
Solomon and and a Goppa code, but perhaps a more efficient
decoding algorithm could be found for the Lomont code \#2 than the
product code, making it also useful in practice.

\begin{openprob}[Decoding Lomont Codes]
Find an efficient algorithm to decode the codes from theorems
\ref{t:p1codes} and \ref{t:genus1degree0Codes}. A good start is
looking at current Goppa code decoding algorithms, say from
Pretzel's book \cite{Pretzel}.
\end{openprob}

\begin{openprob}[Genus 2 classification]
Similar to the codes in theorems \ref{t:p1codes} and
\ref{t:genus1degree0Codes}, classify codes resulting from surfaces
ruled over curves of genus 2.
\end{openprob}

\subsection*{Curves of Garcia-Stichtenoth}\label{s:GSCurves}
% Chris Lomont PhD Thesis 2003
% Curves of Garcia-Stichtenoth

Since good families of codes require long codes, and the curves
explicitly described in this section give Goppa codes reaching the
Drinfeld-Vladut bound, it might be useful to study ruled surfaces
over these curves. In particular, since equations are given below
for the family of curves, and the number of rational points and
genus is known for each curve, computer algorithms could be
developed to study properties of the resulting surface codes.

In 1995, A. Garcia and H. Stichtenoth gave explicit construction
of these two families of curves (GS curves) meeting the TVZ bound
(theorem \ref{t:TVZbound}) \cite{GS1},\cite{GS2}. In the language
of function fields, the G-S curves form towers
${F_1,F_2,F_3,\dots,F_n,\dots}$ of Artin-Schreier extensions of
the rational function field $\field{q^2}(x_1)$. The first family
is given by $$\begin{aligned}
F_1 & = & \field{q^2}(x_1)\\
F_{n+1} & = & F_n(x_{n+1}), \quad n\geq 1\\
x^q_{n+1}+x_{n+1} & = &\frac{x^q_n}{x^{q-1}_n+1}, \quad n\geq 0
\end{aligned}$$
For the first family, the genus  $g(F_n)$ of $F_n$ is shown to be
$$g(F_n)=\left\{
\begin{array}{r@{\quad:\quad}l}
(q^{i+1}-1)(q^{i-1}-1), & \mbox{$n=2i+1$ odd}\\
(q^i-1)^2, & \mbox{$n=2i$ even}\\
\end{array}
\right.
$$
and the smooth curve corresponding to the function field $F_n$ has
at least $(q-1)q^n$ \field{q} rational points. Thus
$\limsup_{n\rightarrow \infty}N_n/g_n(F_n)=q-1$, the
Drinfeld-Vladut bound \cite{DV83}.

The second family is given by
$$
\begin{array}{rcl}
F_1 & = & \field{q^2}(x_1)\\
F_{n+1} & = & F_n(z_{n+1}), \quad n\geq 1\\
z_{n+1}^q + z_{n+1} & = & x^{q+1}_n\\
x_nx_{n-1} & = & z_n \quad n\geq 2
\end{array}
$$
The genus and number of points is known for the second family as
well, and it too reaches the Drinfeld-Vladut bound.

\begin{comment}
 In the second family, the genus
$g(F_n)$ of $F_n$ is shown to be
$$g(F_n)=\left\{
\begin{array}{r@{\quad:\quad}l}
(q^{(n+1)/2}-1)^2, & \mbox{n odd}\\
(q^{n/2}-1)(q^{(n+2)/2}-1), & \mbox{n even}\\
\end{array}
\right.
$$
Thus $g(F_n)\approx q^{n+1}$. Also for the second family, let
$N_n$ be the number of \field{q^2} rational points on the smooth
model of the curve $C$ with function field $F_n$. In odd
characteristic, for $n\geq 3$
$$N_n = (q^2-1)q^{n-1}+2q$$
In even characteristic, for $n\geq 5$
$$N_n = (q^2-1)q^{n-1}+2q^2$$
\end{comment}

It would be worthwhile to study general rank 2 vector bundles over
these curves, in order to construct the corresponding surfaces and
then codes. For decomposable vector bundles, the answers were
given above in section \ref{s:decomposable}, and are product
codes, so cannot beat the best current codes. But for more general
vector bundles programs could be written to investigate the
parameters, perhaps resulting in some good surface codes.

\todo{something wrong with formulas above - fix!}

% end of file

\section{Other Surface Types}
There are codes obtainable from any other surface types, and
perhaps by applying theorem \ref{maintheorem} some code parameters
could be deduced. Some study of the classification of surfaces
would be needed.

\section{Higher Dimensional Varieties} The theory of higher
dimensional varieties is not as well developed as that for curves
and surfaces, but there would be codes here too, although I expect
the problem of determining the code parameters much more difficult
than in the curve or surface case. From a suggestion of Kenji
Matsuki, it seems that studying codes using the code definition
\ref{d:codedef} applied to $\proj^n$ may be solvable.

\begin{openprob}[$\proj^n$ codes]
From the definition \ref{d:codedef} applied to $\proj^n$, find
code parameters and a decoding algorithm.
\end{openprob}

\section{New Projective Curves}\label{s:NewCurveConclusion}
The search in chapter \ref{c:newcurves} for explicit curves
meeting known bounds on number of rational points also has several
possible extensions. Previously known results are \cite{Moreno95}
and \cite{Serre}. The bounds in this paper could possibly be
strengthened by by analyzing the singularities in more detail,
resulting in more known \field{q} rational points on the smooth
models of the curves. Also all genus $\leq 5$ curves from the
degree 6 polynomials were not identified. More work could be done
to compute exact parameters for these curves.

Since computing power grows quickly, the range of curves searched
should be easily done by the time this is in print. For example:
\begin{openprob}[Find more curves]
Extend the techniques of chapter \ref{c:newcurves} on new curves
to search larger spaces, like all degree 7 curves or extend the
coefficient field to \field{4} or larger.
\end{openprob}

% end of file

%\chapter{Inverse Square Root} \todo{include paper??}
%\chapter{Proof K33 nonplanar} \todo{include paper??}
%\chapter{Fast gradient fills} \todo{include paper??}
%\newpage

  % Summary and/or conclusions are optional but often used.
  % The summary and/or conclusions often are the last
  % major division(s) of the text.
  % Reference: PU 19.
% \include{summary}

  % Recommendations are optional.
  % You may include recommendations as a major division if the
  % subject matter and research so prescribe.
  % Reference: PU 19.
% \include{recommend}

  % Bibliography is required if you consulted any outside
  % references.
  % Reference: PU 19.
%\include{bibliography} % todo - include
% end matter for thesis - bibliography, etc

\blankpage

% end of document

  % Appendices are optional.
  % Appendices are not necessarily a part of every thesis.  An appendix is
  % used for supplementary illustrative material, original data, computer
  % print-outs, and other material that is not necessarily appropriate for
  % including in the text of the thesis.
  % Reference: PU 20.
\appendix
%\include{appendix}  % todo - include
% proof of needed results
\chapter*{APPENDIX - SHEAF RESULTS}
\addtocounter{chapter}{1} \addcontentsline{toc}{chapter}{APPENDIX
- SHEAF RESULTS}
\label{S:sheaf} % use this ref to refer to this appendix
\addtocounter{section}{1} Here are collected a few results needed
in this thesis that probably appear in the literature, but I was
unable to find them.

We use the following definition of the binomial coefficient:

\begin{defn}{Binomial Coefficient}\cite[5.1]{GrahamKnuth}
$$\binom{r}{k}= \left \{
\begin{array}{ll}
\frac{r(r-1)(r-2)\dots(r-k+1)}{k!} & k\geq 0 \\
0 & k<0
\end{array}\right . $$

for $r\in\mathbb{R}$ and $k$ an integer.
\end{defn}
Then,
\begin{lemma}\label{l:binlemma}
$$\sum\limits_{j=0}^{n}\binom{j+m}{m}=\binom{m+n+1}{m+1}$$
\end{lemma}

\begin{proof}
$\sum\limits_{j=0}^n\binom{j}{m}=\binom{n+1}{m+1}$, called
``\emph{upper summation}'' from \cite[5.10]{GrahamKnuth}.

Then
$$
\begin{array}{rcl}
\sum\limits_{j=0}^{n}\binom{j+m}{m}&=&\sum\limits_{k=m}^{m+n}\binom{k}{m}\\
&=&\sum\limits_{k=0}^{m+n}\binom{k}{m}-\sum\limits_{k=0}^{m-1}\binom{k}{m}\\
&=&\binom{m+n+1}{m+1}-\binom{m}{m+1}\\
&=&\binom{m+n+1}{m+1}
\end{array}
$$
where the last term in the second to last line is $0$ by
definition.

\end{proof}

\begin{thrm}
\label{tensordeg} For vector bundles \sheaf{E} and \sheaf{F} over
a curve $C$
$$\deg
(\sheaf{E}\otimes\sheaf{F})=\rank\sheaf{F}\deg\sheaf{E}+\rank\sheaf{E}\deg\sheaf{F}$$
\end{thrm}

\begin{proof}
Let \sheaf{E} have rank $r$, and \sheaf{F} have rank $s$. Then by
\cite[Appendix A, C4]{Hartshorne} the Chern polynomials are
$$c_t(\sheaf{E})=\prod\limits_{i=1}^r(1+a_it)$$
and
$$c_t(\sheaf{F})=\prod\limits_{j=1}^s(1+b_jt),$$
where the $a_i$ and $b_j$ are formal symbols. From \cite[Appendix
A, C5]{Hartshorne} or \cite[Remark 3.2.3 (b),(c)]{Fulton} we then
use
$$c_t(\sheaf{E}\otimes\sheaf{F})=\prod_{i,j}(1+(a_i+b_j)t)$$
and the Chern polynomial of an exterior power gives
$$
\begin{array}{rcl}
c_t(\wedge^{r+s}(\sheaf{E}\otimes\sheaf{F})) & = & 1+t\sum\limits_{i,j}(a_i+b_j) \\
 & = & 1+t(s\sum\limits_i a_i + r\sum\limits_j b_j)
 \end{array}
$$

Since the degree of the sheaf in the Riemann-Roch theorem
\cite[Appendix A, Example 4.1.1]{Hartshorne} comes solely from the
$s\sum a_i+r\sum b_j$ term, and since
$c_t(\wedge^r\sheaf{E})=1+t\sum a_i$, etc.,  we get the result.
\end{proof}

Next we compute the rank and degree of $S^n(\sheaf{E})$ where
\sheaf{E} is a vector bundle of rank $r$ and degree $d$ over a
curve $C$.

\begin{thrm}\label{symmetricrank} \textbf{If} \sheaf{E} is a rank $r$ degree $d$
vector bundle over a curve $C$, \textbf{then}
$$\rank S^n(\sheaf{E})=\binom{n+r-1}{r-1}$$
\end{thrm}
\begin{proof}
Since rank is local, and \sheaf{E} is locally free, let $R$ be a
local ring, and compute rank $S^n(M)$, where
$M=\bigoplus\limits_{i=0}^r R$. $S(M)\cong R[x_1,x_2,\dots,x_r]$
\cite[Cor A2.3c]{Eisenbud}. Then $\rank S^n(M)=\{\mbox{\# of
monomials of degree $n$}\}=\binom{n+r-1}{r-1}$.
\end{proof}

The degree is harder to compute. If $\rank\sheaf{E}=1$, then
\sheaf{E} is a line bundle, and $\deg S^n(\sheaf{E}) = n
\deg\sheaf{E}$. So assume $\rank\sheaf{E}\geq 2$. First, applying
the splitting principle \cite[Remark 3.2.3]{Fulton} to \sheaf{E},
$$\sheaf{E}=\sheaf{E}_r\supseteq\sheaf{E}_{r-1}\supseteq\dots\supseteq\sheaf{E}_1\supseteq\sheaf{E}_0=0$$
with $\rank\sheaf{E}_i=i$, line bundles
$\sheaf{L}_i=\sheaf{E}_i/\sheaf{E}_{i-1}$, for $i=1,2,\dots,r$.
Then applying \cite[II, exercise 5.16(c)]{Hartshorne} to
$0\rightarrow\sheaf{E}_{r-1}\rightarrow\sheaf{E}\rightarrow\sheaf{L}_r\rightarrow
0$, we have the filtration
$$S^n(\sheaf{E})=F^0\supseteq F^1\supseteq \dots \supseteq F^n\supseteq F^{n+1}=0$$
and sequences
$$
\begin{array}{cclclclllcc} % oh yeah!
0&\rightarrow& F^1&\rightarrow& S^n(\sheaf{E})&\rightarrow&
S^0(\sheaf{E}_{r-1})&\otimes& S^n(\sheaf{L}_r)&\rightarrow& 0\\
0&\rightarrow& F^2&\rightarrow& F^1&\rightarrow&
S^1(\sheaf{E}_{r-1})&\otimes& S^{n-1}(\sheaf{L}_r)&\rightarrow& 0\\
&&&&&\vdots\\
0&\rightarrow& F^n&\rightarrow& F^{n-1}&\rightarrow&
S^{n-1}(\sheaf{E}_{r-1})&\otimes& S^1(\sheaf{L}_r)&\rightarrow& 0\\
0&\rightarrow& F^{n+1}&\rightarrow& F^n&\rightarrow&
S^n(\sheaf{E}_{r-1})&\otimes& S^0(\sheaf{L}_r)&\rightarrow& 0
\end{array}
$$
From \cite[II, exercise 5.16 (d)]{Hartshorne} degree is additive
across exact sequences, so we get that
\begin{equation}\label{A:sndecomp}
\deg
S^n(\sheaf{E})=\sum\limits_{k=0}^n\deg(S^k(\sheaf{E}_{r-1})\otimes
S^{n-k}(\sheaf{L}_r))\end{equation} Hartshorne states
\cite[proving prop. 2.3]{Hartshorne71} without proof the following
theorem (and it was not clear in his proof what restrictions were
on the bundles other than being over a curve):
\begin{thrm}\label{symmetricdeg}
Let \sheaf{E} be a vector bundle over a curve, of rank $r$ and
degree $d$. \textbf{Then}
$$\deg
S^n(\sheaf{E})=\frac{dn}{r}\binom{n+r-1}{r-1}$$
\end{thrm}

\begin{note}
\emph{Kenji Matsuki noted that this proof can be done in the much
simpler case $\sheaf{E}=\oplus\sheaf{L}_i$, a sum of line bundles,
by using the splitting principle. I will leave the sequences
above, since they may provide tools in some cases to compute the
dimension $k=h^0(C,S^a(\sheaf{E})\otimes\sheaf{O}_C(bP_0))$ for
general vector bundles over curves, as needed in theorem
\ref{ruledtheorem}.}
\end{note}

\begin{proof}
We induct using the above formulas. For rank 1 it is true from
equation \ref{A:sndecomp} above. Assume $\deg\sheaf{E}_{r-1}=d_1$
and $\deg\sheaf{L}_r=d_2$, giving $d=d_1+d_2$. Then
$$
\begin{array}{rcl}
\deg S^n(\sheaf{E}) &=&
\sum\limits_{k=0}^{n}\deg(S^k(\sheaf{E}_{r-1})\otimes
S^{n-k}(\sheaf{L}_r))\\
 &=&\sum \binom{k+r-2}{r-2}(n-k)d_2+\frac{k d_1}{r-1}\binom{k+r-2}{r-2}\\
 &=&\sum (-d_2+\frac{d_1}{r-1})(\frac{k+r-1}{r-1})(r-1)\binom{k+r-2}{r-2} +
 [n d_2-(r-1)(-d_2+\frac{d_1}{r-1})]\binom{k+r-2}{r-2}\\
 &=&\sum (d_1-(r-1)d_2)\binom{k+r-1}{r-1}
 +\sum[n d_2-(d_1-(r-1)d_2)]\binom{k+r-2}{r-2}
\end{array}
$$
where to get to the last line we used the identity
$\frac{k+r-1}{r-1}\binom{k+r-2}{r-2}=\binom{k+r-1}{r-1}$, $r\neq
1$ \cite[5.5]{GrahamKnuth}. Applying lemma \ref{l:binlemma} twice
the sum becomes
$$
\begin{array}{rcl}
&=&
[d_1-(r-1)d_2]\binom{n+r}{r}+[nd_2-(d_1-(r-1)d_2)]\binom{n+r-1}{r-1}
\\
&=&
[\frac{n+r}{r}(d_1-(r-1)d_2)+nd_2-d_1+(r-1)d_2]\binom{n+r-1}{r-1}
\\
&=& [\frac{nd_1}{r}+\frac{nd_2}{r}]\binom{n+r-1}{r-1} \\
&=& \frac{nd}{r}\binom{n+r-1}{r-1}
\end{array}
$$
Where we again used \cite[5.5]{GrahamKnuth} in the form
$\binom{n+r}{r}=\frac{n+r}{r}\binom{n+r-1}{r-1}$. This completes
the proof.
\end{proof}
% end of file

  % Notes and footnotes are optional.
  % Reference: PU 20.
  % I have not implemented this yet.  Mark Senn 2002.06.03
% \include{notes}

  % Vita is required only in a doctoral dissertation.
  % Reference: PU 20.
% Vita for Chris C. Lomont
%
%  vita.tex     February 11, 2002     Mark Senn
%
%  This is the vita for a simple, example thesis.
%
\begin{vita}
\emph{``A little nonsense now and then is relished by the wisest
men"} - Willy Wonka

Chris Lomont was born December 12, 1968, in Fort Wayne, Indiana,
where he spent his youth ``experimenting", which consisted of
explosions and fire. 18 years later, after navigating the public
school system, much to his surprise he found himself studying
math, physics, and computer science in college. He graduated with
bachelor's degrees in all three from Oral Roberts University in
1991, and, finding that working as a waiter was not very exciting,
he started a master's degree in math and/or physics. After mixing
in several programming jobs (including video game programming), he
obtained a master's degree in mathematics from Purdue University,
Fort Wayne campus, in 1996. Certain that school was more fun than
working, he drove to Purdue University, West Lafayette, to pursue
a PhD in superstring theory, which is why seven years later this
thesis is about error correcting codes. While chasing this
lifelong dream of becoming a ``scientist", or, equivalently,
avoiding the private sector, he met the beautiful and brilliant
Melissa Jo Wilson. They were (will be) married on April 26th,
2003. This document is the final hurdle to obtaining his PhD,
which hopefully will be awarded in May 2003. With no more degrees
to chase (on the horizon!), he entered the job market, and the
rest is history.

\end{vita}
 % todo - include


\begin{thebibliography}{}
%%%%%%%%%%%%%%%%%%%%%%%%%%%%%%%%%%%%%%%%%
\vspace{12pt}
\bibitem{AlmHirs92} J. d'Almeida, A. Hirschowitz, \emph{Quelques
plongements non-sp\'{e}ciaux de surfaces rationelles}, Math. Z.
\textbf{211}, pp. 479-483, (1992).

\bibitem{AEJ} J.K. Arason, R. Elman, and B. Jacob, \emph{On indecomposable vector bundles}, Comm. Algebra \textbf{20}, no. 5, pp. 1323-1351, (1992).

\bibitem{Atiyah57} M.F. Atiyah, \emph{Vector bundles over an elliptic curve}, Proc. London Math. Soc., \textbf{7}, pp. 414-452, (1957).

\bibitem{BallCop97} Edoardo Ballico and Marc Coppens, \emph{Very
ample line bundles on blown-up projective varieties}, Bull. Belg.
Math Soc. \textbf{4}, pp. 437-447, (1997).

\bibitem{CRSS98} A. R. Calderbank, E. M. Rains, P. W. Shor
and N. J. A. Sloane, \emph{Quantum Error Correction Via Codes Over
GF(4)}, IEEE Trans. Information Theory, 44 (1998), pp. 1369-1387;
http://xxx.lanl.gov/quant-ph/9608006.

\bibitem{CRSS97} A. R. Calderbank, E. M. Rains, P. W. Shor, and N.
J. A. Sloane, \emph{Quantum error correction and orthogonal
geometry}, Phys. Rev. Lett. \textbf{78}, pp. 405-408, 1997;
http://xxx.lanl.gov/abs/quant-ph/9605005.

\bibitem{Cop95} Marc Coppens, \emph{Embeddings of general
blowing-ups at points}, J. Reine Agnew. Math., \textbf{469} pp.
179-198, (1995).

\bibitem{Cop02a} Marc Coppens, \emph{Very ample linear systems on
blowings-up at general points of projective spaces}, Canad. Math.
Bull., \textbf{45}(3), pp. 349-354, (2002).

\bibitem{Cop02b} Marc Coppens, \emph{Very ample linear systems on
blowings-up at general points of smooth projective varieties},
Pacific Journal of Mathematics, \todo{find short form}
\textbf{202}, No 2, pp. 313-327, (2002).

\bibitem{DV83} V. G. Drinfeld and S. G. Vl\u{a}du\c{t}, \emph{Number of points of an algebraic
curve}, Func. Anal. And Appl., \textbf{17}, pp. 53-54, (1983).
%\todo{pp 68-69 in the Russian}

\bibitem{Eisenbud} D. Eisenbud, \emph{Commutative Algebra with a view towards Algebraic Geometry}, GTM\#150, Springer-Verlag, (1995).

\bibitem{Fulton} W. Fulton, \emph{Intersection Theory}, Second Edition, Springer, (1998).

\bibitem{GS1} A. Garcia and H. Stichtenoth, \emph{A tower of Artin-Schreier extensions of function fields attaining the Drinfeld-Vl\u{a}du\c{t} bound},
Invent. Math., vol 121, pp. 211-222, 1995.

\bibitem{GS2} A. Garcia and H. Stichtenoth, \emph{Asymptotically good towers of function fields over finite fields},
C. R. Acar. Sci. Paris S\'{e}r. I Math., \textbf{322}, pp.
1067-1070, (1996).

\bibitem{Golay49} M. J. E. Golay \emph{Notes on digital coding}, Proc. I.R.E,
\textbf{37}, p. 657, (1949).

\bibitem{Goppa} V.D. Goppa, \emph{Codes on algebraic curves}, Sov.
Math.-Dokl., Vol. 24, pp. 170-172, (1981).

\bibitem{GrahamKnuth} R. Graham, D. Knuth, O. Patashnik, \emph{Concrete Mathematics, A Foundation for Computer
Science}, Addison-Wesley, (1989, 1994).

\bibitem{Grothendieck} A. Grothendieck, \emph{Sur la classification
des fibres holomorphes sur la sph\'{e}re de Riemann}, Amer. J.
Math. \textbf{79}, pp. 121-38, (1957).

\bibitem{Hache} G. Hach\'{e} and D. Le Brigand, \emph{Effective construction
of algebraic-geometric codes}, IEEE Transactions on Information
Theory, vol. 41, pp. 1615-1628, (Nov 1995).

\bibitem{Hamming50} R. W. Hamming, \emph{Error detecting and error
correcting codes}, Bell Syst. Tech. J., \textbf{26}, 2, April
1950.

\bibitem{Hanson} S. H. Hanson, \emph{The geometry of Deligne-Lusztig
varieties; Higher dimensional AG codes}, Ph. D. Thesis, University
of Aarhus, Department of Mathematical Sciences, University of
Aarhus, DK-800 Aarhus C, Denmark, July 1999.

\bibitem{Hartshorne} R. Hartshorne, \emph{Algebraic Geometry}, GTM\#52, Springer-Verlag, (1977).

\bibitem{Hartshorne70} R. Hartshorne, \emph{Resides and duality;
lecture notes on the work of A. Grothendieck, given at Harvard
63/64}, Berlin-New York, Springer-Verlag, (1966).

\bibitem{Hartshorne71} R. Hartshorne, \emph{Ample bundles on curves}, Nagoya Math. J., \textbf{43}, pp. 73-89, (1971).

\bibitem{HL} E. Howe, K. Lauter, http://arxiv.org/abs/math.NT/0207101

\bibitem{Ihara} Y. Ihara, \emph{Some remarks on the number of rational points of algebraic curves over finite fields}, J. Fac. Sci. Tokyo, 28, pp. 721-724, (1981).

\bibitem{Kant} KANT, http://www.math.tu-berlin.de/algebra

\bibitem{Kuchle94} Oliver K\"{u}chle, \emph{Ample line bundles on blown up surfaces},
http://xxx.lanl.gov/abs/alg-geom/9410011

\bibitem{Lauter} K. Lauter, \emph{Geometric methods for improving the upper bounds of the number of rational points on algebraic curves over finite fields}, Journal of Algebraic Geometry, 10, pp. 19-36, (2001).

\bibitem{vanLint} J.H. van Lint, \emph{Introduction to Coding Theory}, GTM\#86, Springer-Verlag, (1982).

\bibitem{Lomont} C. Lomont,  \emph{Yet more projective curves over \field{2}}, to appear in Exp. Math.

\bibitem{LomontSource} www.math.purdue.edu/$\sim$clomont/Math/Papers/2000/PolySolver.zip

\bibitem{DNACode} D. A. Mac D\'{o}naill, \emph{A parity code interpretation of nucleotide
alphabet composition}, Chemical Communications, pp. 2062 - 2063,
(2002).

\bibitem{Maple} MAPLE, www.maplesoft.com

\bibitem{Moreno} C. Moreno, \emph{Algebraic Curves over Finite Fields}, Cambridge Tracts in Mathematics, Cambridge University Press, (1991).

\bibitem{Moreno95} O. Moreno, D. Zinoviev, V. Zinoviev, \emph{On several new
projective curves over $F_2$ of Genus 3, 4, and 5}, IEEE
Transactions on Information Theory, vol. 41, pp. 1643-1648, (Nov
1995).

\bibitem{Peek85} J. B. Hans Peek, \emph{Communications aspects of the compact disc audio
system}, IEEE Communications Magazine, Vol 23, No 2, pp. 7-15 (Feb
1985).

\bibitem{Poonen} Bjorn Poonen, \emph{Bertini theorems over finite fields}, preprint.

\bibitem{Posner68} E. C. Posner, \emph{Combinatorial structures in
planetary reconaissance}, in Error Correcting Codes, edited by H.
B. Mann, pp. 15-46, New York-London-Sydney-Toronto: Wiley (1968).

\bibitem{Pretzel} O. Pretzel, \emph{Codes and Algebraic Curves}, Oxford Science Publications, Clarendon Press, (1998).

\bibitem{Ragot} J. Ragot, http://pauillac.inria.fr/algo/seminars/sem97-98/ragot.html

\bibitem{ReedSolomon60} I. S. Reed and G. Solomon, \emph{Polynomial codes over certain finite
fields}, SIAM J. Appl. Math., Vol 8, No 2, pp. 300-304 (June
1960).

\bibitem{Roman} S. Roman, \emph{Coding Theory and Information Theory}, GTM\#134, Springer-Verlag, (1992).

\bibitem{Schw61} R. L. E. Schwarzenberger, \emph{Vector bundles on the projective
plane}, Proc. London Math. Soc., (3) 11, pp. 623-40, (1961).

\bibitem{Serre55} J.-P. Serre, \emph{Faisceaux alg\'{e}briques
coh\'{e}rents}, Ann. of Math., \textbf{61}, pp. 1-42 (1955).

\bibitem{Serre} J.-P. Serre, \emph{Nombres de points des courbes
algebriques sur $F_q$}, Seminaire de Theorie des Nombres de
Bordeux, expose 22, pp. 1-8, (1983).

\bibitem{Shannon48} C.E. Shannon, \emph{A mathematical theory of
commumication}, Bell Syst. Tech. J., \textbf{27}, pp. 379-423,
623-656, (1948).

\bibitem{Shor96} P. W. Shor, \emph{Fault-tolerant quantum
computation}, Proc. 35th Ann. Symp. on Fundamentals of Computer
Science(IEEE Press, Los Alamitos, 1996), pp. 56-65;
http://xxx.lanl.giv/abs/quant-ph/9605011.

\bibitem{TV} M. A. Tsfafman and S. G. Vl\u{a}du\c{t},
\emph{Algebraic-Geometric Codes}. Dordrecht/Boston/London: Kluwer,
(1991).

\bibitem{TVZ} M.A. Tsfasman, S.G. Vl\u{a}du\c{t}, and Th. Zink, \emph{On
Goppa codes which are better than the Varshamov-Gilbert bound},
Math Nachr, \textbf{109}, pp. 21-28, (1982).


\bibitem{vdvg} van der Geer/van der Vlugt tables at http://www.science.uva.nl/$\sim$geer/

\bibitem{Walker} R. Walker, \emph{Algebraic Curves}, Dover Publications, N.Y., 1962.


\bibitem{XM} http://www.spectrum.ieee.org/WEBONLY/pressrelease/0701/0701dig.pdf

\bibitem{Xu95} Geng Xu, \emph{Ample line bundles on smooth
surfaces}, J. Reine Agnew. Math., \textbf{469}, pp. 199-209,
(1995).

\end{thebibliography}
\end{document}